\documentclass[10pt,twoside]{amsart}

\usepackage{amssymb}
\usepackage{amsxtra}
\usepackage{graphicx}
\usepackage{dcpic, pictexwd}

\theoremstyle{plain}

\newtheorem{thm}{Theorem}
\newtheorem{theorem}{Theorem}[section]
\newtheorem{prop}[theorem]{Proposition}

\newtheorem{defn}[theorem]{Definition}

\newtheorem{lemma}{Lemma}
\newtheorem{cor}{Corollary}

\newcommand{\ga}{\gamma}

\newcommand{\Q}{\mathbb{Q}}

\newcommand{\C}{\mathbb{C}}
\newcommand{\Z}{\mathbb{Z}}

\newcommand\ra{\rightarrow}

\newcommand{\Si}{\Sigma}

\newcommand{\ben}{\begin{enumerate}}
\newcommand{\een}{\end{enumerate}}

\newcommand{\Oz}{P.\ Ozsv{\'a}th\,}
\newcommand{\Sz}{Z.\ Szab{\'o}\,}

\begin{document}

\title{Rational Blow Downs in Heegaard-Floer Homology}
\author{Lawrence P. Roberts}
\address{Department of Mathematics, Michigan State University,
East Lansing, MI 48824}
\email{lawrence@math.msu.edu}
\thanks {The author was supported in part by NSF grant DMS-0353717 (RTG)}

\maketitle

\noindent In \cite{Fint}, R. Fintushel and R. Stern introduced the rational blow down, a process which could be applied to a smooth four manifold
containing one of a family of codimension $0$, negative definite sub-manifolds. This family of examples were extended in the work of J. Park, \cite{Park}. In each case, the procedure called for removing the negative definite piece and replacing it with a rational homology four ball. Furthermore, if the four manifold had a symplectic structure, and the negative definite piece was appropriately embedded with respect to this structure, M. Symington showed that the resulting manifold would have a symplectic structure. \\
\ \\
\noindent The importance of the procedure derives from the possibility that it preserves some of the Seiberg-Witten\footnote{And had a similar property for Donaldson invariants, but we will not be using them here.} basic classes of the original manifold. It could thus be used to construct new examples of homotopy equivalent but non-diffeomorphic four manifolds. Additionally, through clever topological manipulation, it allowed the computation of the Seiberg-Witten invariants for the result of $p$-log transforms on certain four manifolds. In particular, the computation of the Seiberg-Witten invariants for logarithmic transforms in cusp neighborhoods was then achieved. \\
\ \\
\noindent Recently, further extensions of this procedure have been announced, \cite{Stip}, \cite{Froy}. In \cite{Stip} it is noted that another
family of examples built around Wahl singularities, \cite{Wahl}, could also give rise to ways to alter symplectic four manifolds and that in a general class of manifolds, plumbings of spheres along trees, these Wahl examples, when added to those of Fintushel, Stern, and Park, rounded out all the possibilities for symplectic operations. In \cite{Nico}, \cite{Froy}, \cite{Sala}, a further theoretical generalization has been discovered allowing for the possibility of blowing down in the presence of positive scalar curvature.  \\
\ \\
\noindent This paper studies the effect of similar operations on the Ozsv{\'a}th\, and Szab{\'o}\ four manifold invariants, \cite{Four}. In particular,
we will verify that, as conjectured, all the previously known effects of rational blow downs on Seiberg-Witten or monopole invariants have direct analogs for the Ozsv{\'a}th\, and Szab{\'o}\, four manifold invariants. Pending the results of S. Jabuka and T. Mark, \cite{Mark}, these will
extend the calculations of the four manifold invariants to log transforms, just as the original blow downs did. However, more can be said about when the invariants will be preserved, and, to a topologist, a more hands on account of the structure can be given. In particular, we can describe precisely the effect on the $Spin^{c}$ structures. Furthermore, we
can construct more families of examples which exceed those of \cite{Froy} since the relevant three manifolds will probably be hyperbolic, and provide
other possibilities  beyond those considered in \cite{Stip} (although how to find such families in actual four manifolds is unclear). All of these examples rest upon several of Ozsv{\'a}th\, and Szab{\'o}'s main results in \cite{Plum}, \cite{Doub}, and \cite{Unkn}. \\
\ \\
\noindent To be precise, we will study the effect of generalized blow-ups:

\begin{defn} 
Suppose $Y$ is a rational homology sphere which is the oriented boundary of a rational homology ball $B$ and a negative definite four manifold, $C$. If $X_{B}$ is a four manifold containing $B$ as an embedded codimension $0$ piece, and $\varphi \in\mathrm{Diff}^{+}(Y)$, the $(C, \varphi)$-blow up of $X_{B}$ is $X_{C}= (X_{B} - B)\cup_{\varphi}C$. 
\end{defn}

\noindent We start by collecting some facts from Ozsv{\'a}th\, and Szab{\'o}'s papers about rational homology balls and negative-definite cobordisms. This is followed by a proof of the main result and some simplifications of its statement. We then verify that the Seiberg-Witten results for basic classes can be established in this context. Having this theory in hand we go about analyzing all the various examples that were previously known and
drawing out some of the conclusions. Finally, we introduce and analyze some new examples.

\section{Heegaard-Floer results for rational homology balls}

\noindent Let $B$ be a rational homology four ball with boundary $\partial B =Y$, a rational homology sphere. We start by citing a result (Prop 9.9) from \cite{AbsG}:

\begin{prop}\cite{AbsG}
The map $F^{\infty}_{B,\mathfrak{u}}$ is an isomorphism for all $Spin^{c}$ structures on $B$.
\end{prop} 
\ \\
\noindent We will need the following standard

\begin{defn} \cite{Lens} \label{defn:Lspc}
A torsion $Spin^{c}$ structure $\mathfrak{s}$ on  $Y$ is an $L$-structure if $HF_{\mathrm{red}}(Y,\mathfrak{s}) \cong 0$. $Y$ is an $L$-space if it is a rational homology sphere and all its $Spin^{c}$ structures are $L$-structures.
\end{defn}

\noindent {\bf Note:} From \cite{Lens} the set of $L$-spaces is closed under connect sum and has the property that if $K \subset Y$ is a framed knot and $Y_{\infty}$ and $Y_{0}$ are both $L$-spaces then so is $Y_{1}$, where $Y_{i}$ is the result of $i$-surgery on $K$. All three manifolds with elliptic geomoetry are $L$-spaces, but there are others. \\
\ \\
\noindent If we employ the long exact sequence for the flavors of Heegaard-Floer homology we obtain

\begin{lemma}
The maps $F^{-}_{B,\mathfrak{u}}$ and $F^{+}_{B, \mathfrak{u}}$ are isomorphisms onto $HF^{\pm}(Y, \mathfrak{t})/HF_{red}$. If $\mathfrak{s} = \mathfrak{u}|_{Y}$ is an $L$-structure these maps are isomorphisms onto $HF^{\pm}(Y, \mathfrak{t})$. In particular, if $Y$ is an $L$-space, then these maps are isomorphisms for every $\mathfrak{u}$.
\end{lemma}

\noindent {\bf Proof:} The diagram below includes the vanishing of the connecting homomorphism for $S^{3}$. Since the middle map is an isomorphism, and the lower sequence splits, we immediately see that $F^{-}_{B, \mathfrak{u}}$ and $F^{+}_{B, \mathfrak{u}}$ are non-trivial. \\
$$
\begindc{\commdiag}[10]
\obj(0,0)[I1]{$\cdots$}
\obj(10,0)[Ad]{$HF^{-}(S^{3})$}
\obj(20,0)[Ad1]{$HF^{\infty}(S^{3})$}
\obj(30,0)[Ad2]{$HF^{+}(S^{3})$}
\obj(40,0)[I2]{$\cdots$}
\obj(0,10)[I3]{$\cdots$}
\obj(10,10)[Bl]{$HF^{-}(Y,\mathfrak{s})$}
\obj(20,10)[Bl1]{$HF^{\infty}(Y,\mathfrak{s})$}
\obj(30,10)[Bl2]{$HF^{+}(Y,\mathfrak{s})$}
\obj(40,10)[I4]{$\cdots$}
\mor{I1}{Ad}{$0$}
\mor{Ad}{Ad1}{$i_{\ast}$}
\mor{Ad1}{Ad2}{$\pi_{\ast}$}
\mor{Ad2}{I2}{$0$}
\mor{I3}{Bl}{$\delta_{\ast}$}
\mor{Bl}{Bl1}{$i_{\ast}$}
\mor{Bl1}{Bl2}{$\pi_{\ast}$}
\mor{Bl2}{I4}{$\delta_{\ast}$}
\mor{Ad}{Bl}{$F^{-}_{(B,\mathfrak{u})}$}[-1,0]
\mor{Ad1}{Bl1}{$F^{\infty}_{(B,\mathfrak{u})}$}[-1,0]
\mor{Ad2}{Bl2}{$F^{+}_{(B,\mathfrak{u})}$}[-1,0]
\enddc
$$
\ \\
\noindent The isomorphism property follows from the fact that $d (Y,\mathfrak{u}|_{Y}) = 0$ since $(B, \mathfrak{u})$ induces a rational homology, $Spin^{c}$-cobordism to the standard structure on the three sphere. Furthermore, the grading shift induced by any $Spin^{c}$ structure on $B$ will be $0$. Thus for any $\mathfrak{s}$ found by restriction from a $Spin^{c}$ structure on $B$, the map $F^{+}$ will be an isomorphism, if it is non-trivial, onto $HF^{+}(Y, \mathfrak{s})/HF_{red}$. The analogous statement for $F^{-}$ then follows from the exact sequence above and the structure theorem for the tower portion for torsion $Spin^{c}$ structures. When $\mathfrak{s}$ is an $L$-structure, only the tower remains, and the map is an isomorphism. When $Y$ is also an $L$-space, all the cobordism maps will be isomorphisms. $\Diamond$\\
 \\
\noindent This result should be compared with the use of the reducible solution to the Seiberg-Witten equations in the papers by R. Fintushel and R. Stern and Jongil Park. In fact, a negative definite four manifold with $L$-space boundary will have a similar property. The map $F_{W - B^{4},\mathfrak{u}}^{\infty}$ will be an isomorphism for each $\mathfrak{u}$. \\
 \\
\noindent The invariant $d(Y, \mathfrak{s})$ plays an important role below, so it is useful to recall a theorem by Owens and Strle which strengthens the implications of rational homology cobordism:

\begin{lemma} \cite{Owen}
Let $Y$ be a rational homology sphere and suppose that $Y$ bounds a four manifold $X$ which is a rational homology ball. Then $|H_{1}(Y;\,\Z)| = t^{2}$ for some $t \in \mathbb{N}$. There is a $Spin^{c}$ structure on $Y$, denoted $\mathfrak{t}_{0}$, so that $d(Y, \mathfrak{t}_{0} + \beta) =0$ for all $\beta \in \mathrm{Im}(H^{2}(X;\,\Z) \ra H^{2}(Y;\,\Z))$. Furthermore, if $Y$ is a homology lens space, in particular with cyclic furst homology, then given an ordering of the $Spin^{c}$ structures $\{\mathfrak{t}_{0}, \ldots, \mathfrak{t}_{t^{2} -1}\}$ reflecting the action of $H^{1}$, there is a $Spin^{c}$ structure $\mathfrak{t}_{j}$ so that $d(Y, \mathfrak{t}_{j + kt}) = 0$ for $k = 0, \ldots, t - 1$.
\end{lemma}

\noindent It is also worthwhile to note at this point that the for a $\varphi \in \mathrm{Diff}^{+}(Y)$ we will have 

\begin{enumerate}
\item Let $\mathfrak{s} \in Spin^{c}(Y)$, then $HF^{+}_{red}(Y, \varphi^{\ast}\mathfrak{s}) \cong HF^{+}_{red}(Y, \mathfrak{s})$ as the homologies are diffeomorphism invariants as $\Z[U]$-modules. If $\mathfrak{s}$ is an $L$-structure on $Y$, so is $\varphi^{\ast} \mathfrak{s}$,
\item Furthermore, $\varphi$ preserves grading, so $d(Y,\varphi^{\ast}\mathfrak{s}) = d(Y, \mathfrak{s})$
\end{enumerate}

\section{Heegaard-Floer Homology for Certain other Negative Definite Cobordisms}\label{app:plum}
\ \\
\noindent Consider first the general situation where the $L$-space, $Y$, bounds two negative definite four manifolds $W_{1}$ and $W_{2}$. Suppose $W_{1}$ is a codimension $0$ piece of $X_{W_{1}}$ and let $X_{W_{2}}$ be the result of any replacement of $W_{1}$ by $W_{2}$. Suppose further that there are two $Spin^{c}$ structures $\mathfrak{u}_{i}$ on $W_{i}$ which restrict after gluing to the same $Spin^{c}$ structure on $Y$. Let $S(W,\mathfrak{u}) = c_{1}(\mathfrak{u})^{2} + b_{2}(W)$. Assuming everything has trivial rational first homology, we have\\

\begin{lemma}
Suppose $S(W_{1},\mathfrak{u}_{1}) \leq S(W_{2}, \mathfrak{u}_{2})$ and $\Phi_{X_{W_{1}},\mathfrak{u}_{1}\#\mathfrak{u}}(U^{t}) \neq 0$ then $\Phi_{X_{W_{2}}, \mathfrak{u}_{2}\#\mathfrak{u}}(U^{t + n}) \neq 0$ where $8n = S(W_{2}, \mathfrak{u}_{1}) - S(W_{1}, \mathfrak{u}_{2})$. 
\end{lemma}
\ \\
\noindent The proof is a straightforward matching of maps on the various pieces. Since the image of these two maps lie in the same $Spin^{c}$ structure in an $L$-space, the difference is in fact divisible by $8$, so $n$ is an integer. This is the analog of the theorem in \cite{Froy} when applied to $Y$ with positive scalar curvature. The crucial question is when can we match the degree shifts exactly. In Heegaard-Floer homology several powerful results already exist which do just that: they characterize certain negative definite manifolds which are ``sharp'' for the degree inequalities. For them more can be said, the interchange can preserve basic classes, and the restriction to $L$-spaces can be weakened. This is shown in the next section, after a review of the properties of some of these sharp manifolds.\\
\ \\
\noindent Let $W$ be a negative definite cobordism $S^{3} \ra Y$ built by adding $2$-handles to $S^{3} \times I$ along {\em unknots}, which may well
be linked. Usually, we will assume that we have a handlebody presentation for $W$. We let $I$ be the intersection form for $W$ and $Char(I)$ be the characteristic covectors for $I$, i.e. those linear forms $K$ such that $\langle K, v \rangle \equiv v^{t}Iv$ for all $v\in \mathrm{Span}_{\Z}\{h_{i}\}$. We denote $v^{t}Iv$ by $v \cdot v$. We think of $K$ as an element of $H^{2}(W;\,\Z)$, and define $PD[v]$ using Poincare-Lefshetz duality on $W$. We then let $\mathbb{H}^{+}(I) = Hom(Char(I),\mathcal{T}^{+}_{0})$ where we require that the maps have finite support in $Char(I)$ and satisfy the relationship: if $2n = \langle K, v \rangle + v \cdot v $ and $n \geq 0$ then $U^{n}\phi(K + 2PD[v]) = \phi(K)$, but if $n \leq0$ then $\phi(K + 2PD[v]) = U^{-n}\phi(K)$. By acting through $\mathcal{T}^{+}_{0}$ we think of this set as a $\Z[U]$-module.\\
 \\
\noindent $\mathbb{H}^{+}(I) \cong \bigoplus _{\mathfrak{t} \in  Spin^{c}(\partial W)} \mathbb{H}^{+}(I, \mathfrak{t})$ where we can consider $Spin^{c}(\partial W)$ to be the orbits in $H^{2}(W;\,\Z)$ under the action of $2H^{2}(W, \partial W;\,\Z)$, a lattice spanned by the vectors $PD[h_{i}]$, for all the homology classes determined by $h_{i}$, the two handles. The sets $\mathbb{H}^{+}(W,\mathfrak{t})$ are those maps supported on the orbit for $\mathfrak{t}$. \\
\ \\
\noindent Then $W$ induces a map 

$$
T_{W} : HF^{+}(-Y) \ra \mathbb{H}^{+}(I)
$$ 
by
$$
T_{W}(\xi) = \phi \hspace{.5in} \mathrm{where} \hspace{.5in} \phi(K) = F^{+}_{W,\mathfrak{u}_{K}}(\xi)
$$
where $F^{+}_{W,\mathfrak{u}_{K}}$ is the cobordism map induced by the $Spin^{c}$ structure corresponding to $K$. This satisfies the conditions
to be in $Char(I)$ by the generalized adjunction property of \cite{Symp} applied to the spheres which represent the homology classes $[h_{i}]$.  Furthermore, $T_{W}$ maps $HF^{+}(Y,\mathfrak{s})$ to $\mathbb{H}^{+}(I, \mathfrak{s})$. \\
\ \\
There is a grading on $\mathbb{H}^{+}(I)$ given by specifying that a map $\phi$ is homogeneous of degree $d$ if for each characteristic vector $K$ with $\phi(K) \neq 0$, $\phi(K)$ is homogeneous and satisfies:
$$
\mathrm{deg}(\phi(K)) - \left(\frac{K^{2} + b_{2}(W)}{4} \right) = d
$$
\ \\
\noindent where $K^{2}$ is computed using $I^{-1}$. 

\begin{defn} $W$ will be called {\em sleek} when $\partial W$ is a rational homology sphere, $T_{W}$ is a grading preserving isomorphism, and for any $Spin^{c}$ structure $\mathfrak{s}$ on $\partial W$ we have
$$
d(Y,\mathfrak{s}) = \max_{\{K \in \mathrm{Char}_{\mathfrak{t}}(I)\}} \frac{K^{2} + b_{2}(W)}{4}
$$
\end{defn}

\noindent By the fundamental result in \cite{Plum}, this definition is not meaningless.

\subsection{Known Results for sleek Manifolds}

\noindent We first give a synopsis of the result in \cite{Plum}. Let $G$ be a tree with an ordering of its vertices and a map $m : V(G) \ra \Z$, defined on the vertices. Let $d(v)$ be the valence of $v$. We can construct a symmetric matrix $I=(a_{ij})$ by setting $a_{ii} = m(v_{i})$ and $a_{ij} = 1$ if $v_{i}$ and $v_{j}$ are joined by an edge, $0$ otherwise. This is the intersection matrix for $H_{2}(W(G);\,\Z)$, where $W(G)$ is the simply connected four manifold obtained by considering $G$ to define a plumbing of disc bundles over 2-spheres with Euler numbers $m(v_{i})$. $G$ is said to be negative definite when I is negative definite. We call a vertex of $G$ a ``bad'' vertex if $d(v) > -m(v)$. \\
\ \\
\noindent The main theorem of \cite{Plum} states that
\begin{thm}
Let $G$ be a negative-definite graph with at most one bad vertex, then $W(G)$ is a sleek negative definite four manifold.
\end{thm}
\ \\
\noindent There is also a generalization to the case of $G$ with two bad vertices. The manifold is no longer sleek, but one does recover that $T_{W}$ induces an isomorphism $\mathbb{H}^{+}(G,\mathfrak{s}) \cong HF^{+}_{even}(-Y, \mathfrak{s})$, where the even grading is with reference to the canonical $\Z/2\Z$-grading.\\
\ \\
\noindent Also found in \cite{Plum} is an algorithm for finding the rank of $\mathrm{Ker}\ U$ on $HF^{+}(-Y,\mathfrak{s})$ that identifies this set with the set of full paths of characteristic vectors ending in condition (2), according to
 
\begin{defn} A full path of vectors, $K_{0}$, $K_{1}$, $\ldots$,  $K_{n}$, is one where 
\ben
\item[] $\bullet$ $K_{0}$ satisfies
$$
m(v_{i}) + 2\leq \langle K_{0}, v_{i} \rangle \leq -m(v_{i})
$$
\item[] $\bullet$ $K_{i+1}$ is $K_{i} + 2PD[v_{j}]$ where $\langle K_{i}, v_{j} \rangle = - m(v_{j})$. When this occurs we will write $K_{i+1} \sim K_{i}$.
\item[] $\bullet$ $K_{n}$ satisfies either condition (1): there is a vertex $v$ such that $\langle K_{n}, v \rangle > - m(v)$, or condition (2): for all vertices $m(v_{i}) \leq \langle K_{n}, v_{i} \rangle \leq -m(v_{i})-2$.
\een
\end{defn}

\noindent It should be noted that all the vectors in a full path induce the same grading change in their respective cobordism maps. \\
\ \\
\noindent There is another source of sleek negative definite four manifolds, which we will use later. In \cite{Doub}, \cite{Unkn} \Oz and \Sz\,
study the homology of the branched double cover of non-split links in $S^{3}$. Let $\Si(L)$ be the double cover of the link $L$. In particular,
they prove

\begin{prop}
For an alternating non-split link $L$, there is a sleek negative-definite cobordism, $X_{L}$, from $S^{3}$ to $\Si(L)$. 
\end{prop}

\noindent We will review the construction of $X_{L}$ when we need it later in the paper. When $L$ is an alternating knot, $\Si(L)$ is an $L$-space and $X_{L}$ can be found so that $I$ is its Goeritz form. This can be generalized to quasi-alternating links, \cite{Doub}. It should be mentioned that 
the algorithm described above, based on full paths, does not apply to this case. \\ 
\ \\
\noindent Finally, it should be noted that in \cite{Rust} another extension is given for $W(G)$ which is negative semi-definite with one dimensional kernel. It should be possible to combine this with results of \cite{AbsG} to obtain a similar theorem to the one in the next section, interchanging the plumbing piece with a rational homology $S^{2} \times D^{2}$, but the statement is not as clean, and the author knows of no non-trivial examples.

\section{The Effect of Certain Blow-ups}

\begin{thm} \label{thm:main}
Let $X_{B} \subset W$ be a closed, oriented, smooth four manifold with $b_{2}^{+}(X_{B}) > 1$ equipped with a $Spin^{c}$-structure, $\mathfrak{u}$. Suppose $X_{B}$ contains, as a codimension $0$ submanifold, a rational homology ball $B$ whose boundary $Y = \partial B$ satisfies

\ben
\item $Y$ is a rational homology sphere
\item If $\mathfrak{s} = \mathfrak{u}|_{Y}$ then any element in $HF^{+}_{red}(Y, \mathfrak{s})$ has maximal absolute grading less than $-1$. 
\item $Y$ is the oriented boundary of a sleek negative definite four manifold $W$.
\een 
\noindent then, for any $\varphi \in \mathrm{Diff}^{+}(Y)$ the four manifold $X_{C} = (X_{B}  - B) \cup_{\varphi} W$ possesses a $Spin^{c}$ structure $\mathfrak{u}_{C}$ for which 

$$
\Phi_{(X_{B},\mathfrak{u})} = \pm \Phi_{(X_{C}, \mathfrak{u}_{C})}
$$ 
\end{thm} 

\noindent {\bf Note:} 
\ben
\item We will denote $\mathfrak{u}|_{X_{B} - B}$ by $\mathfrak{u}_{0}$.
\item $\Phi$ is only defined up to sign, and thus the equality at the end is as good as can currently be expected.
\item The statement is that they are equivalent as maps. A short Mayer-Vietoris argument for homology with rational coefficients shows that
 $H_{1}(X_{B} - B;\,\Q) \ra H_{1}(X_{\circ};\,\Q)$ is an isomorphism. Choosing a basis for $H_{1}(X_{B} - B;\,\Z)$ can then be shown to provide
 a basis for $H_{1}(X_{\circ};\,\Z)/\mathrm{Tors}$. 
\item The proof will be given for the case where $T_{W}$ induces an isomorphism with $HF^{+}_{even}(-Y)$, as this subsumes the other cases.
\een
\ \\
\noindent{\bf Proof:} Since $\mathfrak{s}$ extends across a rational homology ball, $d(Y,\mathfrak{s}) = 0$. Let $\mathfrak{t} = \varphi^{\ast}\mathfrak{s}$, then $d(Y,\mathfrak{t})= 0$. By our assumptions on $Y$, $HF^{+}(-Y, \mathfrak{t}) \cong\mathcal{T}^{+}_{0} \oplus HF^{+}_{red}(- Y, \mathfrak{t})$, since $d(-Y,\mathfrak{s}) = -d(Y,\mathfrak{s})$, and thus $HF^{+}_{even}(-Y, \mathfrak{t})$ includes the image of $HF^{\infty}(- Y, \mathfrak{t})$. \\
\ \\
\noindent Since $d(Y, \mathfrak{t}) = 0$
$$
d(Y, \mathfrak{t}) = \max_{\{K \in \mathrm{Char}_{\mathfrak{t}}(G)\}} \frac{K^{2} + b_{2}(W)}{4}
$$
there is some $K$ where $K^{2} = - b_{2}(W) = - b_{2}^{-}(W)$ and the $Spin^{c}$ structure $\mathfrak{u}_{K}$ restricts to $\mathfrak{t}$ on $Y$. \\
 \\ 
\noindent Let $W = A \cup B$ where the gluing occurs along a rational homology sphere. Let $c \in H^{2}(W;\,\Z)$ be a characteristic element. R. Fintushel and R. Stern, \cite{Fint}, argue that this occurs if and only if the restriction of  $c$ to each of $A$ and $B$ is characteristic. Thus, there is a characteristic element on $X_{C}$ found by gluing $K$ to the restriction of $\mathfrak{u}$, by way of $\varphi$. \\
 \\
\noindent If we consider the map $F^{+}_{W, \mathfrak{u}_{K}}:\mathcal{T}^{+}_{0} \oplus HF^{+}_{red}(-Y, \mathfrak{t}) \ra\mathcal{T}^{+}_{0}$, we see that for $K$ the grading shift is $0$. Furthermore, this occurs for any $K$ with the properties so far described. \\
 \\
\noindent We now use the map $T_{W}$ and apply it to $\xi$ the generator of the $0$-graded part of $\mathcal{T}^{+}_{0}$ in $HF^{+}(-Y,\mathfrak{t})$. We can find this generator by taking $U^{n} \xi_{2n}$ where $\xi_{2n}$ is the generator of the $2n$-graded part, for large enough $n$. Since $T^{+}$ is an isomorphism of graded $\Z[U]$-modules we see that:

$$
\mathrm{deg}\,\phi_{\xi}(K) = \frac{K^{2} + b_{2}(W)}{4}
$$
\ \\
\noindent when $\phi_{\xi}(K) \neq 0$. Suppose $\phi_{\xi}(K) = 0$ for all $K$ with $K^{2} = - b_{2}(W)$. For every other $K'$, restricting to $\mathfrak{t}$ on $Y$, the right hand side is $< 0$. Thus, $F^{+}_{W, \mathfrak{u}_{K'}}(\xi) = 0$ since the non-zero elements of $\mathcal{T}^{+}_{0}$ have grading $ \geq 0$. Thus $\phi_{\xi}\equiv 0$ as a map in $\mathbb{H}^{+}(I, \mathfrak{t})$. This contradicts that $T_{W}$ is an isomorphism. \\
 \\
\noindent Thus, there must be a $K$, with $K^{2} = -b_{2}(W)$ so that $F^{+}_{W,\mathfrak{u}_{K}}(\xi) = \xi' \neq 0$ where $\xi'$ is in the degree $0$ part of $\mathcal{T}^{+}_{0}$. Furthermore, $U^{n} F^{+}_{W, \mathfrak{u}_{K}}(\xi_{2n}) = \xi'$. $\xi'$ must be primitive, for if it equalled $k \nu$, then $\frac{1}{k}\phi$ would also define a map in $\mathbb{H}^{+}(I, \mathfrak{t})$. This would correspond to an element in the $0$-graded part of the tower in $HF^{+}(- Y, \mathfrak{t})$ such that multiplied by $k$ we obtain $\xi$, a contradiction.  The map $F^{+}_{W, \mathfrak{u}_{K}}$ will then be an isomorphism of towers. However, it may also be non-trivial on an element of $HF^{+}_{red}(- Y, \mathfrak{s})$ \\
 \\
\noindent Dualizing this map we can calculate $F^{-}_{W,\mathfrak{s}_{K}}$ as a map from $HF^{-}(S^{3}) \ra HF^{-}(Y,\mathfrak{t})$. We employ the duality pairing in \cite{Smoo}. Let $\Psi_{(-2)}$ be the generator of the $\mathcal{T}^{-}_{-2}$ tower as a $\Z[U]$-module, and let $\Theta_{(-2)}$ be the same for the tower in $HF^{-}(S^{3})$. Let $\Xi$ be any element in $HF^{-}_{red}(Y, \mathfrak{t})$. Then the duality pairing states:

$$
\langle F^{+}_{W, \mathfrak{u}_{K}}(\xi), \Theta_{(-2)} \rangle = \langle \xi, F^{-}_{W, \mathfrak{u}_{K}}(\Theta_{(-2)}) \rangle
$$
\ \\
\noindent but $\langle \xi', \Theta_{(-2)} \rangle = \pm 1$, so $F^{-}_{W(G), \mathfrak{u}_{K}}(\Theta_{(-2)})=$  $\pm \Psi_{(-2)} + D$ where $D$ is some element in the reduced homology. However, since $\Xi$ has grading $ < -2$, as the connecting homomorphism lowers grading by $1$, and $ F^{-}_{W, \mathfrak{u}_{K}}$ preserves grading, we have that $F^{-}_{W, \mathfrak{u}_{K}}(\Theta_{(-2)})$ does not have a portion in $HF^{-}_{red}(Y, \mathfrak{t})$. The image of $U^{n} \cdot \Theta_{(-2)}$ lies in the tower, $\mathcal{T}^{-}_{-2}$, by $\Z[U]$-equivariance. Thus the composition $\varphi_{\ast} \circ F^{-}_{W, \mathfrak{s}_{K}}$ will be a grading preserving isomorphism onto the tower part of $HF^{-}(Y, \mathfrak{s})$. \\
 \\
\noindent On the other hand, we have already seen that, under our assumptions, the image of $F^{-}_{B,\mathfrak{u}|_{B}}$ is a grading preserving isomorphism onto the $\mathcal{T}^{-}_{-2}$ portion of $HF^{-}(Y, \mathfrak{s})$. If we let $\mathfrak{u}_{0}$ be the restriction of $\mathfrak{u}$ to $W -B$, then from appendix \ref{app:four} we know that

$$
\Phi_{X_{B}, \mathfrak{u}} = \Phi_{X_{B} - B, \mathfrak{u}_{0}} \circ  F^{-}_{B, \mathfrak{u}|_{B}}(\Theta_{(-2)})
$$
\ \\
\noindent and 

$$
\Phi_{X_{C}, \mathfrak{u'}\#_{\varphi}\mathfrak{s}_{K}} = \Phi_{X_{B}-B, \mathfrak{u}_{0}} \circ  \varphi_{\ast} \circ F^{-}_{W, \mathfrak{s}_{K}}(\Theta_{(-2)})
$$
\ \\
\noindent We have shown that both of the right hand sides reduce to

$$
\pm \Phi_{X_{B} - B, \mathfrak{u}_{0}}(\Psi_{(-2)})
$$
\ \\
\noindent where $\Psi_{(-2)}$ generates the $\Z[U]$-module $HF^{-}(Y,\mathfrak{s})$. $\Diamond$ \\
 \\
\noindent {\bf Note:} By the degree formula for a homogeneous $\phi$ we have
$$
0 \leq \mathrm{deg}\,\phi(K) = d + \left(\frac{K^{2} + b_{2}(G)}{4}\right)
$$
\ \\
\noindent For $K$ restricting to $\mathfrak{s}$ as above, we know that $K^{2} + b_{2}(W) \leq 0$, thus $d \geq 0$. Each element
of $HF^{+}_{red}(-Y, \mathfrak{s})$ induces a homogeneous $\phi$ of the same grading under the grading preserving
isomorphism $T^{+}$. Thus, the minimal grading in $HF^{+}_{red}(- Y, \mathfrak{s})$ is $0$. By duality, this implies
that the maximal grading in $HF^{-}_{red}(Y, \mathfrak{s})$ is $-2$. The assumption in the theorem corresponds to requiring
that it be $< -2$. \\
\ \\
\noindent The four manifold invariant, $\Phi$, was extended to four manifolds with $b^{+}_{\,2} = 1$ in \cite{Symp}. In this case, we need
a choice of a line $\mathcal{L} \subset H_{2}(X,\, \Q)$ where $w \in \mathcal{L}$ has $w \cdot w = 0$. One
then splits $X = X_{1} \#_{N} X_{2}$ along a three manifold $N$ with the property that $H_{2}(N,\,\Q) \stackrel{i_{\ast}}{\ra} \mathcal{L}$ $\subset H_{2}(X,\,\Q)$. It is shown in \cite{Symp} that only the choice of $\mathcal{L}$, not $N$, affects the invariant $\Phi_{X, \mathfrak{u}, \mathcal{L}}$. The argument in theorem \ref{thm:main} extends immediately to this setting as long as $N$ can be chosen to avoid $B$ (or $C$). However, a Mayer-Vietoris
argument shows that $H_{2}(X_{B};\,\Q) \cong H_{2}(X_{B} - B;\,\Q)$ so this will always be possible. We have

\begin{cor}
In the setting of Theorem \ref{thm:main}, but with $b_{\,2}^{+}(X_{B}) = 1$ and a choice of a line, $\mathcal{L} \subset H_{2}(X_{B} - B;\,\Q)$, of square zero homology classes, for each $\mathfrak{u}$ on $X_{B}$ there is at least one $\mathfrak{u}_{C}$ on $X_{C, \varphi}$ extending $\mathfrak{u}_{0}$ and having the property that, for every $\varphi \in \mathrm{Diff}^{+}(Y)$, 
$$
 \Phi_{X_{C, \varphi}, \mathfrak{u}_{C}, \mathcal{L}} = \pm \Phi_{X_{B}, \mathfrak{u}, \mathcal{L}}
$$
\end{cor}  
\ \\
\noindent We can simplify the conditions of the theorem in certain circumstances. For instance, as a consequence of \cite{Symp}, if $G$ is a negative definite plumbing graph along a tree that additionally has no bad vertices, then $Y$ is an $L$-space. Using this fact, in \cite{Symp}, \Oz and \Sz prove

\begin{prop}
Let $Y$ be a Seifert fibered rational homology sphere with Seifert invariants $(b, \beta_{1}/\alpha_{1},$ $\,\ldots,\,\beta_{n}/\alpha_{n})$ with $\alpha_{i} > 2$ and $0 < \beta_{i} < \alpha_{i}$ and $(\alpha_{i}, \beta_{i}) = 1$. Then $Y$ bounds a plumbing of sphere bundles defined by a weighted graph, $G$, that is star-like with central node having weight $b$ and whose $i^{th}$ ray is labelled with the ``Hirzebruch-Jung'' fractional expansion of $\beta_{i}/\alpha_{i}$ . If $b \leq -n$ then $G$ is negative-definite and $Y$ is an $L$-space.
\end{prop}

\noindent We can then state, as a result of Theorem \ref{thm:main}

\begin{cor}
Let $Y$ be a Seifert fibered rational homology sphere with $ b \leq -n$ that also bounds a rational homology ball $B$. If $B \subset X_{B}$ 
where $X_{B}$ is a closed, oriented four manifold with $b_{\,2}^{+}(X_{B}) \geq 1$ then, for any $\mathfrak{u} \in Spin^{c}(X_{B})$ and for any $\varphi \in \mathrm{Diff}^{+}(Y)$ the four manifold $X_{C} = (X_{B}  - B) \cup_{\varphi} W(G)$ possesses a $Spin^{c}$ structure $\mathfrak{u}_{C}$ for which 

$$
\Phi_{(X_{B},\mathfrak{u})} = \pm \Phi_{(X_{C}, \mathfrak{u}_{C})}
$$ 
\ \\
where, if $b_{\,2}^{+}(X_{B}) = 1$ these are each computed with reference to an appropriate line $\mathcal{L} \subset H_{2}(X_{B} - B;\,\Q)$. 
\end{cor}

\noindent This theorem addresses three manifolds with positive scalar curvature and the negative definite manifolds found from their presentation as Seifert-fibered manifolds, \cite{Lens}. The reader who is interested mostly in examples may want to skim the next few, largely formal, sections.

\section{Ozsv{\'a}th\,- Szab{\'o}\, Simple Type}

\noindent We follow \cite{Mark} in saying 

\begin{defn} 
A $Spin^{c}$ structure $\mathfrak{u}$ on a closed, oriented, smooth four manifold, $X$, is an {\em Ozsv{\'a}th\,- Szab{\'o}\, basic class} if $\Phi_{X, \mathfrak{u}}\not\equiv 0$. $X$ is {\em  Ozsv{\'a}th\,- Szab{\'o}\, Simple Type} if $D(X, \mathfrak{u}) = 0$ for every basic class, see appendix \ref{app:four} for the definition of $D(X,\mathfrak{u})$. 
\end{defn}

\begin{prop}
If $X_{B}$ satisfies our standing  assumptions and $X_{C}$ is simple type, then $X_{B}$ is also simple type.
\end{prop}

\noindent {\bf Proof:} Suppose that $\Phi_{X_{B}, \mathfrak{u}_{B}}\not\equiv 0$ as a map. By our main theorem we can choose an extension of $\mathfrak{u}_{0}$ to $X_{C}$ by a $Spin^{c}$ structure which pulls back to $W$ with characteristic vector $K$ satisfying $K^{2} = - b_{2}(W)$. Call this $Spin^{c}$ structure $\mathfrak{u}_{K}$. \\
 \\
\noindent Furthermore, $\Phi_{X_{C}, \mathfrak{u}_{K}} \not\equiv 0$ and thus $D(X_{C}, \mathfrak{u}_{K}) = 0$. However,

$$
\begin{array}{c} 4 D(X_{C}, \mathfrak{u}_{K}) = c_{1}(\mathfrak{u}_{K})^{2} - (2\,\chi(X_{C}) +  3\,\sigma(X_{C})) \\
 \\
 = c_{1}(\mathfrak{u}_{0})^{2} + K^{2} - (2\,\chi(X_{B}  ) + 2b_{2}(W) + 3\,\sigma(X_{B}) - 3b_{2}(W)) \\
  \\ 
  = 4D(X_{B}, \mathfrak{u}_{B}) + K^{2} - 2b_{2}(W) + 3b_{2}(W) = 4D(X_{B}, \mathfrak{u}_{B})\\
   \\
\end{array}
$$

\noindent hence $D(X_{B}, \mathfrak{u}_{B}) = 0$. Thus $X_{B}$ is simple type. $\Diamond$ \\
 \\
\noindent This argument also shows that $D(X_{C}, \mathfrak{u}_{K}) =$ $D(X_{B}, \mathfrak{u})$ if we extend $\mathfrak{u}_{0}$ as in Theorem \ref{thm:main}. There is a partial converse to this proposition:

\begin{prop}
Suppose that $\mathfrak{u}_{B}$ has $D(X_{B}, \mathfrak{u}_{B}) = 0$, then $\Phi_{X_{C}, \mathfrak{u}} \not\equiv 0$ for an extension of $\mathfrak{u}_{0}$ to $X_{C}$ only if $D(X_{C}, \mathfrak{u}) = 0$
\end{prop}

\noindent {\bf Proof:} Let $\mathfrak{u}_{K}$ be the extension used in our main theorem, with $K^{2} = -b_{2}(W)$ on $W$. Let $\mathfrak{u}_{K'}$ be any other extension of $\mathfrak{u}_{0}$ to $W$. Then $K' - K = 2PD[V]$ for some $V$ in the span of the Poincare duals of the spheres in the cobordism. \\
 \\
\noindent If we calculate $D(X_{C}, \mathfrak{u}_{K})$ and $D(X_{C}, \mathfrak{u}_{K'})$ we find that

\begin{equation} \label{eqn:grrel}
D(X_{C}, \mathfrak{u}_{K'}) - D(X_{C}, \mathfrak{u}_{K}) =\langle K, V \rangle + V\cdot V
\end{equation}

\noindent However, 

$$
\langle K, V \rangle + V\cdot V =\frac{1}{4}\big((K + 2PD[V])^{2} -K^{2} \big) = \frac{1}{4}\big((K + 2PD[V])^{2} + b_{2}(W) \big) \leq 0
$$
\ \\
\noindent The last inequality comes from using $K + 2PD[V]$ in

$$
d(Y, \mathfrak{t}) = \max_{\{K \in \mathrm{Char}_{\mathfrak{t}}(I)\}}\frac{K^{2} + b_{2}(W)}{4}
$$ 
\ \\
\noindent The proof of the previous proposition shows that $D(X_{C}, \mathfrak{u}_{K}) =$ $D(X_{B}, \mathfrak{u}_{B})$, hence if $X_{B}$ is Ozsv{\'a}th\,-Szab{\'o}\, simple type then so is $X_{C}$ for the lifted $Spin^{c}$ structures. $\Diamond$

\section{Other Extensions of the $Spin^{c}$-structure}

\noindent Throughout this section we assume that $\mathfrak{u}|_{Y}$ is an $L$-structure. Results such as we will prove can be adapted 
if all the reduced homology is in sufficiently negative grading relative to $D(X_{B}, \mathfrak{u})$. \\
 \\
\noindent The other $Spin^{c}$ structures on $W$ restricting to $\mathfrak{s}$ on $Y$ occur in the orbit of $\mathfrak{u}_{K}$ under the action of $\mathrm{Ker}(H^{2}(W;\,\Z) \ra H^{2}(Y;\,\Z))$. That is, by the lattice generated by $PD[v_{i}]$ for all $i$. In fact, the $K$ in the theorem can be chosen from the set of $K$ satisfying $v_{i} \cdot v_{i}  \leq\ \langle K, v_{i} \rangle\ \leq - v_{i} \cdot v_{i}$. \\
 \\
\noindent We can extend the maps to the other elements in the orbit of $K$ by the relationship defining $\phi$. Namely, if 

$$
2n = \langle K, v \rangle + v \cdot v
$$
\ \\
\noindent and $n \geq 0$ then $U^{n}\phi(K + 2PD[v]) = \phi(K)$, but if $n \leq 0$ then $\phi(K + 2PD[v]) = U^{-n}\phi(K)$. \\
 \\
\noindent Running back through the implications, we arrive at the fact, also obtained by the adjunction property of \cite{Symp}, that $F^{-}_{W, \mathfrak{u}_{K'}}$ $= U^{n} F^{-}_{W,\mathfrak{u}_{K}}$ for some $n$. It is important that the elements of $H_{2}(W:\Z)$ are represented by spheres. Were $\mathfrak{u}$ not to restrict as an $L$-structure, the action of $U$ could mask potentially deleterious reduced homology; this is the cause for our cautionary assumption. If $K ' - K = 2PD[V]$ where $V \in H_{2}(W;\,\Z)$, we can compute $n$ using equation \ref{eqn:grrel}, the grading shift formula. Since the action of $U$ shifts grading by $-2$, we have 

$$
n = -\frac{\langle K, V \rangle + V\cdot V}{2}
$$
\ \\
\noindent We have shown the following corollary:

\begin{cor} Under the assumptions of Theorem \ref{thm:main}, we have

$$
\Phi_{(X_{B},\mathfrak{u})}(\gamma) = \pm \Phi_{(X_{C},  \mathfrak{u}_{K'})}(U^{\frac{\langle K, V \rangle + V\cdot V}{2}}\cdot \gamma)
$$
\ \\ 
\noindent when $\mathfrak{u}|_{Y}$ is an $L$-structure. Here $\mathfrak{u}_{K'}$ is the structure in $Spin^{c}(X_{C})$ obtained by extending by $K + 2 PD[V]$ instead of $K$.
\end{cor}

\noindent This should be read as applying as long as the $U^{\ast} \cdot \gamma$ still has non-negative degree. Since $\langle K, V \rangle + V \cdot V \leq 0$ this provides an improvement on Theorem \ref{thm:main} when $D(X_{B}, \mathfrak{u}) > 0$. 

\section{Taut Configurations}

\noindent We follow \cite{Fint} in the definition

\begin{defn} 
Let $W \subset X$ be a codimension $0$-submanifold of $X$. We say that $W$ is {\em tautly embedded relative to a structure} $\mathfrak{u} \in Spin^{c}(X)$ if the characteristic class, $K$, of $\mathfrak{u}|_{W(G)}$ satisfies

$$ 
|\langle K, v_{i} \rangle| + v_{i} \cdot v_{i} \leq -2
$$
\ \\ 
\noindent We will say that $W$ is tautly embedded in $X$ if it is tautly embedded relative to the Ozsv{\'a}th\,- Szab{\'o}\, basic classes of $X$.
\end{defn}

\noindent Due to the adjunction inequalities in \cite{Smoo}, \cite{Symp}, this is only a meaningful restriction when $v$ is a sphere of 
negative self-intersection. It is illuminating to consider the case when we have a plumbing of spheres along a tree with corresponding four manifold $W(G)$:

\begin{lemma} 
Suppose that $G$ has $m(v_{i}) < - 1$ for every vertex and that $W(G) \subset X$. That $W(G)$ is tautly embedded relative to $\mathfrak{u}$ implies that there is a Stein structure on $W(G)$ whose induced $Spin^{c}$ structure agrees with $\mathfrak{u}|_{W(G)}$. 
\end{lemma}

\noindent {\bf Proof:} If one draws the Kirby diagram for the tree plumbing, the assumptions allow us to find a Legendrian link diagram for $W(G)$, where $L_{i}$ represents $v_{i}$, the framing on each component is $tb(L_{i}) - 1$, and the rotation number satisfies $rot(L_{i}) = \langle K, [v_{i}] \rangle$. The latter is possible since, by assumption, the $K$ values satisfy the Thurston-Bennequin inequality. In this diagram all the attaching circles are along a horizontal line, with branches from the same vertex nesting, see the solution to exercise 6.3.9(e) of \cite{Four}. By results in \cite{Four}, the first Chern class, $c_{1}(J)$, of the Stein structure given by Legendrian surgery on this link also has $rot(L_{i}) = \langle c_{1}(J), [v_i] \rangle$. $\Diamond$ \\

\begin{lemma}
Each basic class $\mathfrak{u}$ on $X_{B}$ which restricts to an $L$-structure on $\partial B$ can have at most one lift to $X_{C}$ relative to which $W(G)$ will be tautly embedded.
\end{lemma}

\noindent {\bf Proof:} Let $\mathfrak{s} = \mathfrak{u}|_{Y}$ where $Y = \partial B$. Since $\mathfrak{s}$ is an $L$-structure by assumption, there is precisely one full path of characteristic vectors on $W$, restricting to $\varphi^{\ast} \mathfrak{s}$, constructed from the algorithm in \cite{Plum}, and terminating in condition 2 above. If the inequalities 

$$
v_{i}\cdot v_{i} + 2 \leq \langle K, v_{i} \rangle \leq - v_{i} \cdot v_{i} - 2
$$
\ \\
\noindent are satisfied, the full path consists of precisely one vector, $K$, as the algorithm terminates immediately. $\Diamond$ 
\ \\
\noindent This will {\em not} be true outside the tree plumbing setting. \\
\ \\
\begin{prop}
Suppose $W(G)$ is tautly embedded in $X_{C}$, with $\partial W$ an $L$-space, then for each basic class $\mathfrak{u}$ on $X_{B}$, there is a unique lift $\hat{\mathfrak{u}}$ to $X_{C}$ such that

$$
\Phi_{X_{C}, \hat{\mathfrak{u}}} = \pm \Phi_{X_{B}, \mathfrak{u}} 
$$
\end{prop}
\ \\
\noindent On the other hand, we can always choose the lift $K$ to be the original vector in a full path terminating in condition 2. If $HF^{+}_{red} \not\cong 0$, we choose the path with $K^{2} = - b_{2}(W)$. If that path contains more that one vector, then we can find a number of lifts equal to the cardinality of the path. Each of these vectors will have $K^{2} = -b_{2}(W)$ by construction and $$\phi(K + 2PD[v]) = \phi(K)$$ implies that the $\phi$ generated by the lowest degree element in the tower of $HF^{+}(-Y, \mathfrak{s})$ will induce an isomorphism on homologies. These lifts will also preserve the simple type condition. 

\section{Examples}

\subsection{Blow-up/Blow-down Formulas}

\noindent  In this case, $B = B^{4}$ and $G$ is a single vertex with $-1$-multiplicity. If $v$ represents the homology class of the $-1$-framed sphere, and $e$ is the hom-dual cohomology class, then $PD[v] = - e$. We know that $K = (2n + 1) e$. These all lie in the same full path, for which the formulae above provide

$$
\Phi_{X', \mathfrak{s}\#(2n+1)e} (U^{-\frac{n(n+1)}{2}} \cdot \ga) = \Phi_{X,  \mathfrak{s}}(\ga)
$$

\noindent where $\ga = U^{m} \otimes h$. This is the blow-up formula from \cite{Smoo}. \\
 \\
\noindent This comports with the calculation for $D(\hat{\mathfrak{s}}_{n})$ where we have

$$
D(\hat{\mathfrak{s}}_{n}) = D(\mathfrak{s}) - (n^{2} + n)
$$

\noindent If $D(\mathfrak{s}) > 0$ then it is possible for there to be non-zero invariants for $n \neq \pm 1$. 

\subsection{Generalized Rational Blow-Down} 

\noindent We verify the formulas for the generalized rational blow-down as defined by Jongil Park, \cite{Park}, in his extension of the rational blow-down of R. Fintushel and R. Stern, \cite{Fint}. \\
 \\
\noindent Consider the lens space $-L(p^{2}, p\,q - 1)$ where $1 \leq q < p$ and $(p,q) = 1$. A. Casson and J. Harer, \cite{Cass}, showed that this bounds a rational homology ball $B_{p,q}$. Let $\frac{p^2}{p\,q - 1}  = [b_{k},b_{k-1}, \ldots, b_{1}]$ where $b_{i} \geq 2$. Then $-L(p^{2}, p\,q - 1)$ also bounds the four manifold, $C_{p,q}$ found from a linear plumbing of $2$-sphere bundles with Euler numbers $-b_{k}, \ldots, -b_{1}$, as in the diagram

$$
\begindc{\undigraph}[10]
\obj(0,0)[n1]{$- b_{k} $}
\obj(8,0)[n2]{$- b_{k-1}$}
\obj(16,0)[n3]{$-b_{k-2}$}
\obj(30,0)[n4]{$-b_{1}$}
\mor{n1}{n2}{$\ $}
\mor{n2}{n3}{$\ $}
\mor{n3}{n4}{$\cdots$}
\enddc
$$
\ \\
\noindent Suppose $X_{(c)} = C_{p,q} \cup X'$, then the \textit{generalized  rational blow down} of $X_{(c)}$ along $C_{p,q}$ is the four manifold $X_{(b)} = B_{p,q} \cup X'$. The case when $q=1$ is the original \textit{rational blow down} described by R. Fintushel and R. Stern, \cite{Fint}. This 
is uniquely defined as every diffeomorphism of $-L(p^{2}, p\,q - 1)$ extends to the rational homology ball. \\
 \\
\noindent In \cite{Park}, J. Park calculates the various algebraic topological facts about these manifolds. In particular, $\pi_{1}(B_{p,q}) = \Z/p\Z$ and we can choose bases so that the inclusion homomorphism

$$
i^{*} : H^{2}(B_{p,q};\,\Z) \cong \Z/p\Z \ra H^{2}(L(p^{2}, p\,q - 1);\,\Z) \cong \Z/p^{2}\Z
$$
\ \\
is given by $ n \ra np$. This could also be deduced from the Owens-Strle result mentioned earlier. \\
 \\
\noindent As $\partial W(G)$ is a Lens space, the conditions of the theorem apply, verifying that the generalized rational blow down results will hold in the Heegaard-Floer theory. \\
 \\
\noindent {\bf I.} In \cite{Fint} it is shown that there is a copy of $C_{2,1}$, an embedded sphere with self-intersection $-4$, in $E(2)$ for which
$X_{B}$ is diffeomorphic to $3\,{\C}P^2 \# 18\,\overline{{\C}P}^{2}$. The latter divides along $S^{3}$ into two pieces each with $b_{\,2}^{+} > 0$ and
thus $\Phi \equiv 0$ for every $Spin^{c}$ structure. On the other hand $\Phi_{E(2),\mathfrak{u}_{0}} = 1$, \cite{Symp}, for the $Spin^{c}$ structure with trivial first Chern class. For the graph, $G$, with one vertex of multiplicity $-4$, we need only consider $K$ with $K = -4, -2, 0, 2, 4$ when paired with the sphere. For these, $K^{2} + 1 = 0$ only for $K = \pm 2$. However, the basic class for $E(2)$ restricts as $K = 0$. Thus, upon removing $C_{2,1}$ and replacing it with $B_{2,1}$, the restriction of the basic class does not extend. For those classes which do extend, both $\Phi_{3\,{\C}P^2 \# 18\,\overline{{\C}P}^{2}, \mathfrak{u}'} = 0$ and $\Phi_{E(2),\mathfrak{u}'} = 0$, in keeping with the theorem. \\
\ \\
\noindent {\bf II.} We consider the plumbing diagram for $C_{5,3}$ whose boundary is $-L(25, 14)$.

$$
\begindc{\undigraph}[6]
\obj(0,0)[n1]{$-3$}
\obj(10,0)[n2]{$-5$}
\obj(20,0)[n3]{$-2$}
\mor{n1}{n2}{$\ $}
\mor{n2}{n3}{$\ $}
\enddc
$$
\ \\
\noindent For this plumbing 
$$
P^{-1} = -\frac{1}{25}\left( \begin{array}{ccc} 
		14 & 3 & 1 \\
		3 & 6 & 2 \\
		1 & 2 & 9 \\
		\end{array} \right)
$$
\noindent There are $30$ classes with $m(v_{i}) + 2 \leq \langle K, v_{i} \rangle \leq - m(v_{i})$, but only the five classes

$$
\begin{array}{c}
(1,3,0)  \hspace{.25in} (-1,-3,0)  \hspace{.25in} (3,-1,0) \\ 
\ \\
 (-1,1,2) \hspace{.25in} (1, -3, 2)  \\
\end{array}
$$
\ \\
\noindent give rise to characteristic classes with $K^{2} = -3$. The restriction of these to $-L(25,14)$ also extend to the rational ball $B_{5,3}$
since there are $5$ structures which do extend. Two, $(1,3,0)$  and $(-1,-3,0)$, correspond to Stein structures on $W(G)$, while the other three occur in full paths of length $2$.

$$
(3,-1,0) \sim (-3,1,0) \hspace{.25in}
 (-1,1,2) \sim (-1,3,-2) \hspace{.25in}
 (1, -3, 2) \sim (1,-1,-2)
$$
\ \\
\noindent In \cite{Park} There is an embedding of $C_{5,3}$ into $E(3) \# 2\,\overline{{\C}P}^{2}$. However, as there, we know that the basic classes
on $E(3)$ are $\pm PD[f]$, \cite{Mark}, and so the basic classes on $E(3) \# 2\,\overline{{\C}P}^{2}$ are, up to sign, $PD[f] + a_{1} e_{1} + a_{2} e_{2}$ where $a_{i} = \pm 1$. The spheres have the form $s$, $f - 2e_{1} -e_{2}$ and $e_{1} - e_{2}$, so the basic classes applied to these spheres give
$(1, 2a_{1} + a_{2}, -a_{1} + a_{2})$. Choosing $a_{1} = a_{2} = 1$ provides $(1,3,0)$. This and its negative are the only two of the five which can occur, and thus will determine the basic classes on $X_{B}$. It is illuminating to compare this with the calculation in \cite{Park}. \\
\ \\
\noindent {\bf III. Log Transforms} We assume that all our four manifolds are simply connected and are simple type, with $b_{\,2}^{+}(X) > 1$. These can be relaxed somewhat, but merely make the exposition tedious. \\
\ \\
\noindent In \cite{Mark}, S. Jabuka and T. Mark calculate the Ozsv{\'a}th-Szab{\'o} basic classes for elliptic surfaces. They consider the formal generating function

$$
\mathcal{OS}(X) = \sum_{\mathfrak{u} \in Spin^{c}(X)} \Phi_{W, \mathfrak{u}} e^{c_{1}(\mathfrak{u})}
$$
\ \\
\noindent which is an element of $\Z[H^{2}(X;\,\Z)]$ under our simplicity assumptions and prove that  
$\mathcal{OS}(\mathrm{E}(n))$ is supported on the non-zero terms of $(e^{2[T]} - e^{-2[T]})^{n-2}$ where $[T]$ is the Poincar\'e dual of 
the fiber class. As the Ozsv{\'a}th-Szab{\'o} four manifold invariant does not have a canonical sign choice, using the formal generating
functions directly is problematic. However, as a consequence they show that the Ozsv{\'a}th-Szab{\'o} basic are the same as the Seiberg-Witten basic classes.\\
\ \\
We can now follow the argument in \cite{Fint}, see also \cite{Four}, to establish that the Seiberg-Witten invariants of manifolds found by log transform on a fiber class in an elliptic surface equal the Ozsv{\'a}th-Szab{\'o} invariants for each $Spin^{c}$ structure. In \cite{Fint}, it is shown that if $X$ is an irreducible, closed four manifold that is simple type and contains a $0$-framed torus in a cusp neighborhood, there is
an embedding of $C_{p,1}$ in $X \# (p-1)\,\overline{{\C}P}^{2}$, which, when rationally blown down, yields a manifold diffeomorphic to
that obtained by the logarithmic transform of multiplicity $p$ on the torus. We have the following analog of their theorem:

\begin{thm}
Let $X_{B}$ be the result of blowing down this configuration (and thus be the result of a $p$-log transform on the torus), then 
$$
\Phi_{X_{B}, L + s\,F_{p}} = \pm \Phi_{X, L} \hspace{.5in} s = -p + 1,\, - p + 3, \ldots,\, p - 1 
$$ 
where $F_{p}$ is the Poincar{\'e} dual of the multiple fiber in $X_{B}$ obtained from the $p$-log transform. 
\end{thm}
\ \\
\noindent As usual this rests upon the observations that, if $L$ is a basic class for $X$, then $L_{J} = L + \sum J(i)\,e_{i}$ with $J(i) = \pm 1$
are the basic classes for $X \# (p-1)\,\overline{{\C}P}^{2}$ In \cite{Fint} it is shown that $L_{J}|_{X - C}$ has a unique extension to
$X_{B}$ as $L + |J|\,F_{p} \in H^{2}(X_{B};\,\Z)$ where $|J| = \sum J(i)$. This can be verified by the handle arguments below. \\
\ \\
\noindent Together these allow us to compute the invariants for multiple log transforms on distinct fibers of an elliptic surface. In particular,
we know that the Ozsv{\'a}th-Szab{\'o} basic classes for $E(n)_{p,q}$  for $(p,q) = 1$ are the same as the Seiberg-Witten classes, and up to sign, the invariants agree. The specific classes are listed in Theorem 3.3.6 of \cite{Four}.  

\begin{figure}
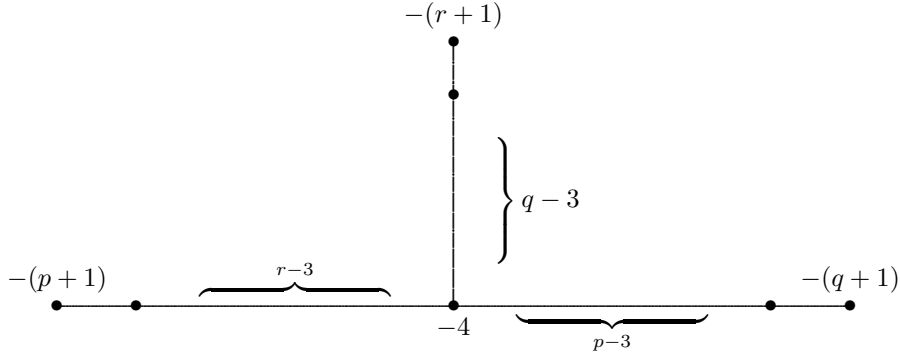

$$
\begindc{\undigraph}[10]
\obj(0,0)[n1]{$- (p + 1) $}
\obj(3,0)[m1]{$\ $}
\obj(15,0)[n2]{$- 4$}[\south]
\obj(27,0)[m3]{$\ $}
\obj(30,0)[n3]{$ - (q + 1)$}
\obj(15,8)[m4]{$\ $}
\obj(15,10)[n4]{$- (r + 1)$}
\mor{n1}{m1}{$\ $}
\mor{m1}{n2}{$\overbrace{\hspace{1in}}^{r - 3}$}
\mor{n2}{m3}{$\underbrace{\hspace{1in}}_{p - 3}$}[\atright, \solidline]
\mor{m3}{n3}{$\ $}
\mor{n2}{m4}{$ \left. \begin{array}{l} \, \\ \, \\ \, \\ \, \\ \end{array} \right\} q - 3$}[\atright, \solidline]
\mor{m4}{n4}{$\ $}
\enddc
$$
\caption{Each unlabelled vertex represents a $-2$-sphere and the braces indicate how many occur along each branch, excluding the one depicted. We assume $p, q, r \geq 2$ and count the number of spheres as $0$ if it is negative.} \label{fig:Wahl}
\end{figure}

\subsection{Wahl-type Plumbing Diagrams}

\noindent In \cite{Stip}, \Sz and A. Stipsicz show that the only plumbing configuration of spheres that can be placed in a symplectic configuration are those of J. Park, considered previously, and a class discovered by J. Wahl, \cite{Wahl}, and depicted in Figure \ref{fig:Wahl}. The latter are also known to bound rational homology balls, \cite{Neum}. When the configuration is symplectic it can be rationally blown down to obtain a symplectic structure on $X_{B}$. That the theorem applies to this configuration follows from the observation above that a tree with no bad vertices represents an $L$-space. We make some observations about the topology of these manifolds and then analyze some of the simplest cases.

\subsubsection{\ }
\noindent In the preceding we relied upon the incidental correspondence between classes with $K^{2} = -|G|$ and those extending to the associated rational homology ball. In the following it will be convenient to have a more reliable way to assert that the classes we find will restrict to structures on $Y$ which extend to the ball. For the examples from \cite{Fint}, \cite{Park}, and \cite{Wahl}, we can find such a criterion by noting that
in each case the rational homology ball is a Mazur manifold. Thus adding a single $2$-handle to $Y$ will give $S^{1} \times S^{2}$, to which we add a $3$-handle and a $4$-handle to obtain the homology ball. A $Spin^{c}$ structure on $Y$ will extend over the rational homology ball as long as it
induces the trivial $Spin^{c}$ structure on $S^{1} \times S^{2}$. We will present $Y$ as $\partial W(G)$, using the relative handlebody calculus in \cite{Four}, and exploit the descriptions in \cite{Cass} and \cite{Neum} to find this new handle. Technically, we are describing $W(G) \cup \overline{B}$, so the argument actually finds those $\mathfrak{s}$ which induce on $-Y$ a $Spin^{c}$ structure that extends over $\overline{B}$. 
Changing orientations shows that $\mathfrak{s}$ will then extend as well.\\
\ \\
\noindent Once we have found this handle, we look at the intersection form for the new framed link and find a primitive element in its kernel. This element is the coefficient vector for a homology class in $H_{2}(W(G) \cup h_{2} ;\,\Z)$ which algebraically intersects none of the classes from the surfaces coming from the link. In this simple setting it must then lie in the boundary, and being primitive, is one of $\pm [S^{2}]$. We can pair $K + a\,h^{\ast}$ with this element, and find those $K$ for which there is an $a$ that makes this pairing $0$. Such $K$'s restrict to $Y$ in a way that will extend over the rational homology ball. As we already know we can limit our search to those $K$ with $m(v_{i}) + 2 \leq \langle K, v_{i} \rangle \leq - m(v_{i})$.

\begin{figure}
\begin{center}
\includegraphics[scale=0.75]{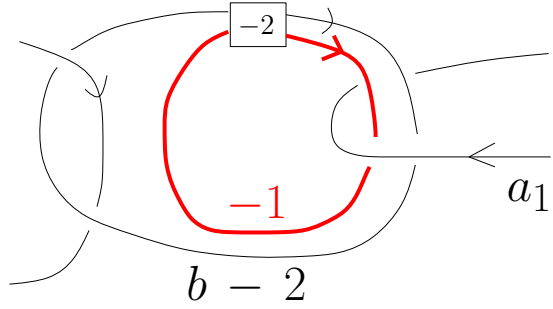}
\end{center}
\caption{The handle addition for constructing a rational homology ball from $C_{p,q}$} 
\label{fig:Cp}
\end{figure}
\ \\
\noindent For Park's examples we can use \cite{Cass} to find the correct handle. For Wahl's examples, we will use \cite{Neum}. In Figure \ref{fig:Cp}
we depict the handle to be added to obtain $B_{p,q}$. To obtain this diagram, we assume that $q < \frac{p}{2}$. This assumption is innocuous since $q + q' = p$ implies that $(p\,q - 1)(p\,q' -1) \equiv 1$ $\mathrm{mod}\ p^{2}$, and thus $-L(p^{2},p\,q - 1) \cong -L(p^{2}, p\,q' - 1)$. The following diagram represents a four manifold with boundary $-L(p^{2}, p\,q - 1)$

$$
\begindc{\undigraph}[7]
\obj(0,5)[n1]{$\frac{q}{p-q}$}[\south]
\obj(5,5)[n2]{$0$}[\north]
\obj(10,5)[n3]{$\frac{p-q}{p}$}[\south]
\obj(5,0)[n4]{$+1$}[\south]
\mor{n1}{n2}{$\ $}
\mor{n2}{n3}{$\ $}
\mor{n2}{n4}{$\ $}
\enddc
$$
\ \\
The assumption on $p$ and $q$ allows us to write 
$$
\begin{array}{c}
\frac{p-q}{p} = [-1,\,-1,\,-3,\,a_{1},\,\ldots] \\
\ \\
\frac{q}{p-q} = [-1,\,-1,\, -\frac{p}{q}]  = [-1,\,-1,\,b,\,\ldots]\\
\end{array}
$$
where $[c, [d]] = c - \frac{1}{[d]}$ iteratively defines the continued fractions. From $a_{1}$ on, the terms are $< -1$. Likewise, $-\frac{p}{q}$
can be written as a continued fraction with entries $<-1$ as $p > q$. Expanding the two fractional surgeries into plumbing chains according
to this recipe, and using the handle in \cite{Cass}, produces the diagram in Figure \ref{fig:Cp}. In the case when $q=1$, $b = -p$, we can check that the diagram above will give the original rational blow-down plumbing diagram for $C_{p,1}$ from \cite{Fint}. In this case, the homology class from the new handle intersects that from the $-(p+2)$-framed handle $-2$ times, and intersects the first $-2$-handle $+1$ times. For the intersection form including this handle 

$$
(-p,\, 1,\, 2 - p,\, 3 - p, \ldots,\, -2,\, -1)
$$
generates the kernel, with the new handle listed first and the other handles in left to right order along the plumbing. If we restrict to taut 
embeddings, then $K(v_{1}) = - p, -p + 2, \ldots, p$ and $K(v_{j}) = 0$. Thus, we need to solve $-p\,a + K(v_{1}) = 0$ for numbers in
the required range. this can only be done if $a = \pm 1$. Hence, we recover a lemma of \cite{Fint}, see also \cite{Four}. \\

\begin{figure}
\begin{center}
\includegraphics[scale=0.75]{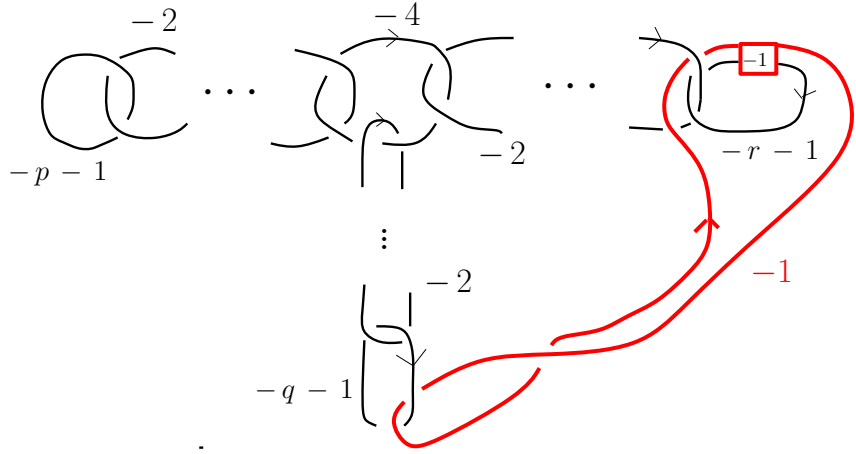}
\end{center}
\caption{Adding the handle in red will give $S^{1} \times S^{2}$. Note there may be more full twists in the handle than depicted; we only need the
intersection data.} 
\label{fig:wahlhandle}
\end{figure}

\noindent Following the diagrams through \cite{Neum} we find the $-1$-framed handle in Figure \ref{fig:wahlhandle} 
that will produce the rational homology ball.  

\subsubsection{$p = 2$, $q = 2$, $r = 2$}
This has plumbing diagram

$$
\begindc{\undigraph}[5]
\obj(0,0)[n1]{$- 3$}[\northwest]
\obj(5,5)[n2]{$- 4$}[\northeast]
\obj(10,0)[n3]{$- 3$}
\obj(5,10)[n4]{$- 3$}
\mor{n1}{n2}{$\ $}
\mor{n2}{n3}{$\ $}
\mor{n2}{n4}{$\ $}
\enddc
$$
\ \\
\noindent We can obtain this configuration in elliptic surfaces with sufficiently non-generic fibers. In particular, blowing up the triple point
of a fiber of type IV in the list in \cite{Kirb} produces the following configuration of proper transforms and exceptional spheres: 

$$
\begindc{\commdiag}[7]
\obj(0,21)[I1]{$-2$}[\southwest]
\obj(10,30)[I2]{$\ $}
\obj(10,21)[J1]{$-2$}[\southeast]
\obj(0,30)[J2]{$\ $}
\obj(5,20)[K1]{$-2$}[\south]
\obj(5,31)[K2]{$\ $}
\mor{I1}{I2}{$\ $}[-1,2]
\mor{J1}{J2}{$\ $}[-1,2]
\mor{K1}{K2}{$\ $}[-1,2]
\obj(13,25)[d1]{$\ $}
\obj(22,25)[d2]{$\ $}
\mor{d1}{d2}{$\# \overline{{\C}P}^{2}$}
\obj(25,25)[A1]{$\ $}
\obj(45,25)[A2]{$-1$}[\east]
\obj(30,20)[A3]{$-3$}[\south]
\obj(30,30)[A4]{$\ $}
\obj(35,20)[A5]{$-3$}[\south]
\obj(35,30)[A6]{$\ $}
\obj(40,20)[A7]{$-3$}[\south]
\obj(40,30)[A8]{$\ $}
\mor{A1}{A2}{$\ $}[-1,2]
\mor{A3}{A4}{$\ $}[-1,2]
\mor{A5}{A6}{$\ $}[-1,2]
\mor{A7}{A8}{$\ $}[-1,2]
\obj(35,18)[e1]{$\ $}
\obj(30,10)[e2]{$\ $}
\mor{e1}{e2}{$\# 3\,\overline{{\C}P}^{2} $}
\obj(7,5)[B1]{$\ $}
\obj(32,5)[B2]{$-4$}[\east]
\obj(17,0)[B3]{$-3$}[\south]
\obj(17,10)[B4]{$\ $}
\obj(22,0)[B5]{$-3$}[\south]
\obj(22,10)[B6]{$\ $}
\obj(27,0)[B7]{$-3$}[\south]
\obj(27,10)[B8]{$\ $}
\obj(9,2)[C3]{$-1$}[\south]
\obj(9,8)[C4]{$\ $}
\obj(11,2)[C5]{$\ $}
\obj(11,8)[C6]{$-1$}
\obj(13,2)[C7]{$-1$}[\south]
\obj(13,8)[C8]{$\ $}
\mor{B1}{B2}{$\ $}[-1,2]
\mor{B3}{B4}{$\ $}[-1,2]
\mor{B5}{B6}{$\ $}[-1,2]
\mor{B7}{B8}{$\ $}[-1,2]
\mor{C3}{C4}{$\ $}[-1,2]
\mor{C5}{C6}{$\ $}[-1,2]
\mor{C7}{C8}{$\ $}[-1,2]
\enddc
$$
\ \\
\noindent If we let $f_{1}$, $f_{2}$, and $f_{3}$
be the three $-2$-spheres then $f = f_{1} + f_{2} + f_{3}$ is the fiber class of the elliptic surface. After blowing up $4$ times, the classes
of the Wahl plumbing in $E(n) \# 4\,\overline{{\C}P}^{2}$ are $e_{1} - e_{2} - e_{3} - e_{4}$, $f_{1} - e_{1}$, $f_{2} - e_{1}$, and $f_{3} - e_{1}$.
When we blow down the Wahl plumbing, we reduce $b_{2}^{-}(E(n) \# 4\,\overline{{\C}P}^{2})$ by $4$.  \\
\ \\
\noindent Since the spheres in the type IV fiber can be made holomorphic in $E(n)$, we have that $PD[f](f_{i}) = 0$. The basic classes
of $E(n)$ are $L_{p} = (2 - n + 2p)PD[f]$ for $p = 0, \ldots, n - 2$, so the basic classes of $E(n) \# 4\,\overline{{\C}P}^{2}$ are then
$$
L_{p} \pm e_{1}^{\ast} \pm e_{2}^{\ast} \pm e_{3}^{\ast} \pm e_{4}^{\ast}
$$
where $e_{i}^{\ast} = -PD[e_{i}]$. These pair with the classes in the Wahl plumbing as $(2,1,1,1)$, $(0,1,1,1)$, $(-2,1,1,1)$, and $(-4,1,1,1)$ and their conjugates. Notice that for $(-4,1,1,1)$, the configuration is {\em not} tautly embedded relative to the corresponding basic classe: $L_{p} - e_{1}^{\ast} + e_{2}^{\ast} + e_{3}^{\ast} + e_{4}^{\ast}$.\\
\ \\
\noindent We can calculate
$$
\begin{array}{ccc}
P =  \left( \begin{array}{cccc} 
		 -4 & 1  & 1  & 1  \\
		 1  & -3 & 0  & 0  \\
		 1  & 0  & -3 & 0  \\
		 1  & 0  & 0  & -3 \\
		\end{array} \right)  & \hspace{0.25in} &

P^{-1} = -\frac{1}{81}\left( \begin{array}{cccc} 
		 27 & 9  & 9  & 9  \\
		 9  & 30 & 3  & 3  \\
		 9  & 3  & 30 & 3  \\
		 9  & 3  & 3  & 30 \\
		\end{array} \right) \\
\end{array}
$$
\ \\
$Y$ will have $|H^{2}(Y;\,\Z)| = 81$, and $9$ of the $Spin^{c}$ structures on $Y$ will extend. We need to consider the $108$ vectors
in $\{ -2, 0,2,4\} \times$ $\{-1,1,3\} \times$ $\{-1,1,3\}\times$ $\{-1,1,3\}$. \\
\ \\
Following our previous prescription, we enhance the intersection form by
$$
P' =  \left( \begin{array}{ccccc} 
		 -4 & 1  & 1  & 1  & 1\\
		 1  & -3 & 0  & 0  & -1\\
		 1  & 0  & -3 & 0  & 1\\
		 1  & 0  & 0  & -3 & 0\\
	   1  & -1 & 1 & 0 & -1 \\
		\end{array} \right)
$$
\ \\
\noindent This matrix has rank $4$ and kernel generated by $(3,-2,4,1,9)$. We need to find those of the $108$ vectors such that
$$3\,a_{1} - 2\,a_{2}+4\,a_{3} + a_{4} \equiv 0\ \mathrm{mod}\ 9$$ We can eliminate many by noticing that this requires 
$a_{2} + a_{3} + a_{4} \equiv 0\ \mathrm{mod}\ 3$. It is then easy to check that the two vectors $(2,1,1,1)$ and $(-2,-1,-1,-1)$, a conjugate pair for which $W(G)$ will be tautly embedded, satisfy the relationship. Among those for which $W(G)$ is not tautly embedded we can find the remaining $7$: $(-2,3,1,-1)$, $(0,1,3,-1)$, $(4,-1,-1,-1)$, $(-2,1,-1,3)$, $(0,3,-1,1)$, $(-2,-1,3,1)$ and $(0,-1,1,3)$. All have square $-4$. However,

$$
\begin{array}{cccc}
(-2,3,1,-1) \sim (0,-3,1,-1) &
(0,1,3,-1) \sim (2,1,-3,-1) &
(4,-1,-1,-1) \sim (-4,1,1,1) \\
(-2,1,-1,3) \sim (0,1,-1,-3) &
(0,3,-1,1) \sim (2,-3,-1,1) &
(-2,-1,3,1) \sim (0,-1,-3,1) \\
\ & (0,-1,1,3) \sim (2,-1,1,-3) & \ \\
\end{array}
$$
\ \\
\noindent At this point, the $\Z/3\Z$-symmetry of the original plumbing diagram has been lost, due to our breaking the symmetry in the handle addition used to describe the rational homology ball. The boundary of this ball has a $\Z/3\Z$ symmetry, and this handle indicates how, relative to the symmetry, the ball is glued when replacing $W(G)$. In our example, $E(n) \# 4\,\overline{{\C}P}^{2}$, the only basic classes which will extend over the ball are 
$$
\begin{array}{cc}
L_{p} - e_{1}^{\ast} - e_{2}^{\ast} - e_{3}^{\ast} - e_{4}^{\ast} &  L_{p} + e_{1}^{\ast} + e_{2}^{\ast} + e_{3}^{\ast} + e_{4}^{\ast}\\
L_{p} - e_{1}^{\ast} + e_{2}^{\ast} - e_{3}^{\ast} - e_{4}^{\ast} &  L_{p} - e_{1}^{\ast} + e_{2}^{\ast} + e_{3}^{\ast} + e_{4}^{\ast}\\
\end{array}
$$
\ \\
\noindent The last two, however, differ by $2PD[v]$ where $v$ is the $-4$ sphere. Hence they restrict to the complement of $W(G)$ as the same $Spin^{c}$-structure, and extend to $X_{B}$ as the same structure. Thus, $X_{B}$ will have three times the number of basic classes as $E(n)$
determined by how they restrict to $B$ in $E(n) \# 4\,\overline{{\C}P}^{2}$. The manifold $X_{B}$ will be minimal. The configuration is that of blow ups of holomorphic spheres and the exceptional curves, and thus can be made symplectic. Hence $X_{B}$ is symplectic with canonical class being the extension of $L_{0} + e_{1}^{\ast} + e_{2}^{\ast} + e_{3}^{\ast} + e_{4}^{\ast}$, as this is the canonical class for $E(n) \# 4\,\overline{{\C}P}^{2}$.
Through the following steps we see that $X_{B}$ is simply connected.  Since $B$ is a Mazur manifold, it will have cyclic fundamental group equal to its first homology, $\Z/9\Z$. On the other hand, this needs to be a subgroup of the first homology of the boundary, which is $\Z/3\Z \oplus \Z/27\Z$. If we let $m_{1}$ be the meridian of the attaching circle of the $-4$ sphere, and $m_{2}$, $m_{3}$, and $m_{4}$ be those of the other spheres; then the first factor is generated by $m_{2} - m_{3}$ and the second factor is generated by $m_{2}$. When we add the additional handle to define the ball we establish that, in $B$, $m_{2} - m_{3} \equiv m_{1}$. As the former has order $3$ and the latter order $9$, we find that $m_{1}$ has order $3$ in $B$. On the other hand, $3m_{2} = m_{1}$, so $m_{2}$ will have order $9$ and thus generate both the first homology and the fundamental group of $B$. However, in $E(n) - N(f)$, $m_{2}$ is just the intersection of a section with the boundary, and thus contractible. By the Seifert-Van Kampen theorem, $X_{B}$ will be simply connected.  For more meaningful examples of blowing down Wahl type plumbings, see \cite{Mich}. 
 
\subsection{Other Examples and Branched Double Covers}
\ \\
\ \\
\noindent{\bf I.} Returning to the generalized rational blow down, it can be seen that if $-L(p^{2},p\,q-1)$ bounds $B$, a rational homology
ball, then $L(p^{2},p\,q -1)$ bounds $\overline{B}$. So for each of the Lens spaces given above, their reverse orientation can also be used.
For example, $-L(25,4) \cong L(25,21) \cong$ $-L(25,19) \cong L(25,6)$. Thus, $-L(25,6)$ bounds a rational homology ball and the negative definite
manifold 

$$
\begindc{\undigraph}[5]
\obj(0,0)[n1]{$-5$}
\obj(10,0)[n2]{$-2$}
\obj(20,0)[n3]{$-2$}
\obj(30,0)[n4]{$-2$}
\obj(40,0)[n5]{$-2$}
\obj(50,0)[n6]{$-2$}
\mor{n1}{n2}{$\ $}
\mor{n2}{n3}{$\ $}
\mor{n3}{n4}{$\ $}
\mor{n4}{n5}{$\ $}
\mor{n5}{n6}{$\ $}
\enddc
$$
\ \\
\noindent Clearly, the main theorem will now also apply to these cases. In addition, these Lens spaces occur as the double branched covers
of $S^{1}$ with branch locus $10_{3}$ for $L(25,6) \cong L(25,21)$ and $\overline{10_{3}}$ for $L(25,4) \cong L(25,19)$. Several other of the generalized rational blow downs occur as branched covers of slice, alternating knots with fewer fewer than ten crossings. $6_{1}$ is the first, with double branched cover $L(9,2) \cong L(9,5)$; its mirror has branched cover $L(9,7) \cong L(9,4)$. The mirrors of $8_{8}$, $10_{22}$ and $10_{35}$ provide the Lens space boundaries for $C_{5,2}$, $C_{7,2}$ and $C_{7,3}$. Of course, the knots themselves also yield Lens spaces upon finding the 
double branched cover, and these also bound negative definite plumbings and rational homology balls. \\
\ \\
\noindent These considerations suggest how to find more examples. According to \cite{Doub} the branched double cover of an alternating knot bounds a 
negative definite four manifold with the requisite properties. If the knot is also slice, the branched double cover will also bound a rational
homology ball. For instance, the Wahl example above arises as the branched double cover of the alternating slice knot $10_{75}$, with one of its orientations. In this setting the main theorem becomes:

\begin{cor}
Let $K$ be an alternating, slice knot in $S^{3}$. Let $X_{K}$ be the negative definite manifold constructed in \cite{Plum} bounded by $\Si(K)$, and
let $B_{K}$ be the rational homology ball found as the double branched cover of $B^{4}$ over the slice disc. Then, for any $\varphi \in \mathrm{Diff}^{+}(\Si(K))$ the four manifold $X_{C} = (X_{B}  - B) \cup_{\varphi} X_{K}$ possesses a $Spin^{c}$ structure $\mathfrak{u}_{C}$ for which 
$$
\Phi_{(X_{B},\mathfrak{u})} = \pm \Phi_{(X_{C}, \mathfrak{u}_{C})}
$$ 
\end{cor}
\ \\
\noindent{\bf II.} The branched double cover of the knot $8_{9}$ is $L(25,7)$, a lens space to which J. Park's results do not apply. $8_{9}$ is a fully amphichiral, alternating, slice, two-bridge knot. Thus $L(25,7)$ bounds a rational homology ball (the branched cover of $B^{4}$ along the slice disc), a negative definite manifold, and is oriented diffeomorphic to itself with reverse orientation ($L(25,18)$). For this reason we will be lax about orientations in gluing. Furthermore, by the results of Bonahon it has symmetry group $\Z/4\Z$, induced by this diffeomorphism and that from the branched covering. Thus, relative to the orientation from the four manifold we are blowing up/down, there will be only one way to glue the rational homology ball into the manifold. Since this is also a Lens space we already have enough to draw the conclusion of the main theorem. However, it is useful and interesting to delve more into the topology. \\
\ \\
First we record the negative definite plumbing diagram for $-L(25,7)$:
$$
\begindc{\undigraph}[7]
\obj(0,0)[n1]{$-4$}
\obj(10,0)[n2]{$-3$}
\obj(20,0)[n3]{$-2$}
\obj(30,0)[n4]{$-2$}
\mor{n1}{n2}{$\ $}
\mor{n2}{n3}{$\ $}
\mor{n3}{n4}{$\ $}
\enddc
$$
\ \\
\noindent It is straightforward to check, using the methods so far described, which vectors $K$ both initiate full paths and have square $-4$. Written
with their full path they are:
$$
\begin{array}{c}
(4,-1,0,0) \sim (-4,1,0,0) \\ 
\ \\
(2,-1,2,0) \sim (2,1,-2,2) \sim (2,1,0,-2) \\
\ \\
(0,1,0,2) \sim (0,1,2,-2) \sim (0,3,-2,0) \sim (2,-3,0,0) \\
\ \\
(-2,3,0,0) \sim (0,-3,2,0) \sim (0,-1,-2,2) \sim (0,-1,0,-2) \\
\ \\
(-2,-1,0,2) \sim (-2,-1,2,-2) \sim (-2,1,-2,0) \\
\end{array}
$$

\begin{figure}
$$
\begindc{\commdiag}[7]
\obj(2,21)[I1]{$-2$}[\southwest]
\obj(12,30)[I2]{$\ $}
\obj(12,21)[J1]{$-2$}[\southeast]
\obj(2,30)[J2]{$\ $}
\obj(7,13)[K1]{$\ $}
\obj(7,31)[K2]{$-2$}[\north]
\obj(0,16)[L1]{$\ $}
\obj(14,16)[L2]{$-4$}
\mor{I1}{I2}{$\ $}[-1,2]
\mor{J1}{J2}{$\ $}[-1,2]
\mor{K1}{K2}{$\ $}[-1,2]
\mor{L1}{L2}{$\ $}[-1,2]
\obj(16,25)[d1]{$\ $}
\obj(24,25)[d2]{$\ $}
\mor{d1}{d2}{$\# \overline{{\C}P}^{2}$}
\obj(25,25)[A1]{$\ $}
\obj(45,25)[A2]{$-1$}[\east]
\obj(30,20)[A3]{$-3$}[\south]
\obj(30,30)[A4]{$\ $}
\obj(35,13)[A5]{$\ $}
\obj(35,30)[A6]{$-3$}[\north]
\obj(40,20)[A7]{$-3$}[\south]
\obj(40,30)[A8]{$\ $}
\obj(28,16)[A9]{$\ $}
\obj(42,16)[A10]{$-4$}
\mor{A9}{A10}{$\ $}[-1,2]
\mor{A1}{A2}{$\ $}[-1,2]
\mor{A3}{A4}{$\ $}[-1,2]
\mor{A5}{A6}{$\ $}[-1,2]
\mor{A7}{A8}{$\ $}[-1,2]
\obj(37,12)[e1]{$\ $}
\obj(30,9)[e2]{$\ $}
\mor{e1}{e2}{$\# 2\,\overline{{\C}P}^{2} $}
\obj(16,6)[B1]{$\ $}
\obj(27,6)[B2]{$-2$}[\east]
\obj(19,4)[B3]{$\ $}[\north]
\obj(19,8)[B4]{$-3$}
\obj(21,0)[B5]{$\ $}
\obj(21,9)[B6]{$-3$}[\north]
\obj(24,3)[B7]{$-2$}[\south]
\obj(24,12)[B8]{$\ $}
\obj(22,10)[B11]{$\ $}
\obj(28,10)[B12]{$-1$}
\obj(26,7)[B13]{$\ $}
\obj(26,13)[B14]{$-5$}
\obj(18,2)[B9]{$-4$}
\obj(25,2)[B10]{$\ $}
\mor{B9}{B10}{$\ $}[-1,2]
\mor{B1}{B2}{$\ $}[-1,2]
\mor{B3}{B4}{$\ $}[-1,2]
\mor{B5}{B6}{$\ $}[-1,2]
\mor{B7}{B8}{$\ $}[-1,2]
\mor{B11}{B12}{$\ $}[-1,2]
\mor{B13}{B14}{$\ $}[-1,2]
\enddc
$$
\caption{\ }
\label{fig:E4L257}
\end{figure}
\ \\
\noindent Since $8_{9}$ is slice, the branched double cover also bounds a rational homology ball, $B$. Later we will present $B$ explicitly. For now
we employ our example. Consider the elliptic surface $E(4)$. It has a section with self-intersection $-4$. We can find a realization of $E(4)$ with 
a type IV fiber, consisting of three $-2$ holomorphic spheres intersecting at a point. Since the section intersects the fiber once homologically, and
all these sub-manifolds are holomorphic, the section intersects only one of the three  $-2$ spheres. We blow up the triple intersection point, and
then twice blow up the intersection of the exceptional spheres with one of the $-3$ spheres not intersecting the fiber. This is schematically described
in Figure \ref{fig:E4L257}. Starting with the $-4$ sphere, and working up and to the right, we encounter a copy of our plumbing diagram. This we will
remove and replace with the rational homology ball. \\
\ \\
\noindent Let $f$ be the homology class of a generic fiber. Then $f = f_{1} + f_{2} + f_{3}$ where $f_{i}$ is the homology class of one of the $-2$ spheres intersecting at the triple point. We choose $f_{1}$ to intersect the section $s$. Since these are holomorphic spheres, any basic class pairs with them to give $0$. Otherwise, the basic classes are $L_{-2} = -2PD[f]$, $0$ and $L_{2} = 2PD[f]$. When we blow up three times the new basic classes are $2PD[f] \pm e_{1}^{\ast} \pm e_{2}^{\ast} \pm e_{3}^{\ast}$ and $e_{1}^{\ast} \pm e_{2}^{\ast} \pm e_{3}^{\ast}$ up to conjugation. The spheres in the plumbing diagram are $s$, $f_{1} - e_{1}$, $e_{1} - e_{2}$, and $e_{2} - e_{3}$. Finally we note that $PD[f_{1}](f_{1}) = -2$, but $PD[f_{i}](f_{1}) = + 1$ for $i = 2, 3$. We can then check by pairing these homology classes with the basic vectors that $L_{2} + e_{1}^{\ast} - e_{2}^{\ast} - e_{3}^{\ast}$, $L_{2} - e_{1}^{\ast} + e_{2}^{\ast} - e_{3}^{\ast}$, and $L_{2}- e_{1}^{\ast} - e_{2}^{\ast} + e_{3}^{\ast}$ restrict to 
$W(G)$ as the three characteristic vectors $K$ in the second full path above. Thus, upon blowing down, we will have an extension with non-trivial
four manifold invariant corresponding to $L_{2}$. By conjugation, there will also be one corresponding to $L_{-2}$ arising from the last full path
listed above. None of the other basic classes on $X_{C}$ will give rise to basic classes on $X_{B}$. In particular the resulting manifold has two basic classes. It is then straightforward to check that $\chi(X_{B}) = 47$ and $\sigma(X_{B}) = -31$. In particular, this is not homotopy equivalent to an elliptic surface or a log transform of one. It's not clear whether this manifold is minimal; however, it is simply connected. We will see below that
a generator of the fundamental group of $B$ is the meridian of the first $-2$ sphere (starting from the $-4$ end). However, from the diagram this will contract along a $-3$ sphere in the complement of $B$. \\

\begin{figure}
\includegraphics[width=3.5in,height=2in]{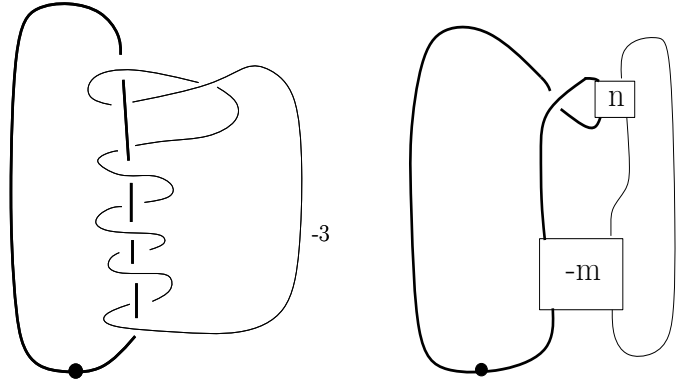}
\caption{The rational homology ball on the left has boundary $-L(25,7)$. It is in the class on the right when $n=1$, $m = 4$. The framing on the two handle on the right is $n - m$.}\label{fig:257fam}
\end{figure}

\begin{figure}
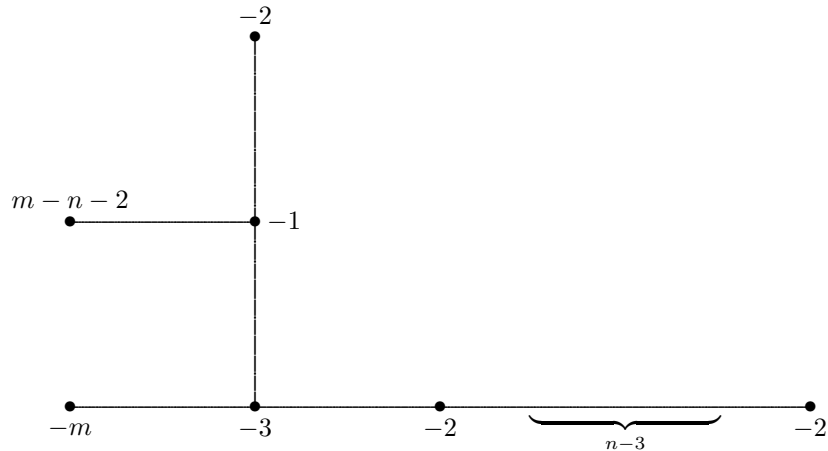

$$
\begindc{\undigraph}[7]
\obj(0,0)[n1]{$-m$}[\south]
\obj(10,0)[n2]{$-3$}[\south]
\obj(20,0)[n3]{$-2$}[\south]
\obj(40,0)[n4]{$-2$}[\south]
\obj(10,10)[n5]{$-1$}[\east]
\obj(0,10)[n6]{$m-n-2$}[\north]
\obj(10,20)[n7]{$-2$}
\mor{n1}{n2}{$\ $}
\mor{n2}{n3}{$\ $}
\mor{n3}{n4}{$\underbrace{\hspace{1in}}_{n - 3}$}[\atright, \solidline]
\mor{n2}{n5}{$\ $}
\mor{n5}{n7}{$\ $}
\mor{n6}{n5}{$\ $}
\enddc
$$
\caption{The boundary of the family of rational homology spheres on the right in Figure \ref{fig:257fam}. The Kirby calculus demonstration
has been relegated to the end of this paper. The meridian of the $m-n-2$ sphere is the generator of the first homology of the rational homology ball. When $n < 3$, we should use a total of $n-1$ spheres along that branch.} \label{fig:257neg}
\end{figure}
\begin{figure}
\includegraphics[width=4in,height=1.5in]{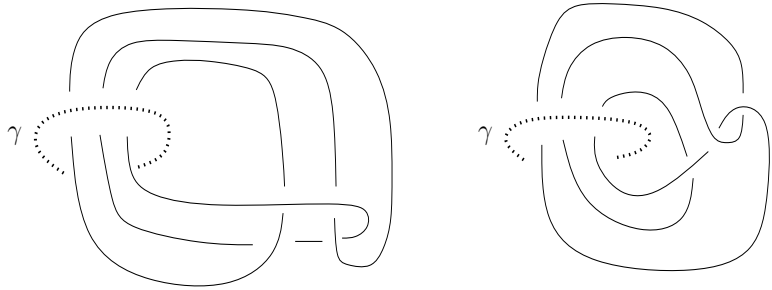}
\caption{The knots in the family are built from the units in this diagram by taking $\Z/n\Z$-cyclic branched covers of $S^{3}$ with axis $\gamma$. In the complement of $\gamma$ we may isotope the knot from the diagram on the left to the one on the right. This is done by moving the outer strand in the left diagram to be the middle strand. Making an equivariant isotopy shows that the knot in the cyclic cover is alternating since the diagram on the right will lift to an alternating projection. On the other hand the diagram on the left clearly lifts to a ribbon projection.} \label{fig:slicefamily}
\end{figure}

\ \\
\noindent We now turn to characterizing the rational homology ball. In fact it sits within an infinite family all of which have interest. In Figure \ref{fig:257fam} we depict the rational homology spheres, whereas in Figure \ref{fig:257neg} we depict plumbing diagrams with the same boundary. There
are several interesting choices of $m$ and $n$. Choosing $m = n+ 3$ gives a family negative definite plumbings with no bad vertices (since $m \geq 1$), the first of which is the lens space $-L(25,7)$. When $m = n + 1$, we can blow down and simplify to finally arrive at the plumbing diagram for $-L(k^{2}, 2 k - 1)$ where $k = 2n + 1$. These are some of Park's examples. When $m = n + 2$, the boundary is no longer irreducible; instead it is the connect sum of $\mathbb{R}P^{3}$ with 

$$
\begindc{\undigraph}[5]
\obj(0,0)[n1]{$-(n+ 2)$}
\obj(10,0)[n2]{$-3$}
\obj(20,0)[n3]{$-2$}
\obj(50,0)[n6]{$-2$}
\mor{n1}{n2}{$\ $}
\mor{n2}{n3}{$\ $}
\mor{n3}{n6}{$\underbrace{\hspace{1in}}_{n - 3}$}[\atright, \solidline]
\enddc
$$
\ \\
\noindent These are still $L$-spaces by the connect sum formula, and they still bound rational homology balls. For the remainder, if we choose $m - n - 2 < -1$ we obtain negative definite plumbings with at most one bad vertex, so the main theorem still applies. \\
\ \\
\noindent  For example, $10_{48}$ has a Seifert-fibered double branched cover which, with one of its orientations, is the boundary of the plumbing diagam in Figure \ref{fig:257neg} with $m = 5$ and $n=2$.\\
\ \\
\noindent{\bf III.} The knot $9_{41}$ is the first non two bridge, non Montesinos, slice, alternating knot. Its double branched cover is neither
a Lens space nor a Seifert fibered space. It is the third in a family of ribbon knots which are all alternating and which we analyze simultaneously. This family appears first in \cite{Ngky} and is used in \cite{Naik} in a study of equivariant concordance. It can be easily described as in Figure \ref{fig:slicefamily} from which it is seen that the first is the unknot, and the second the knot $6_{1}$, analyzed above. This family is also interesting since the branched cover of $9_{41}$ is hyperbolic (verified using SnapPea) and the rest are likely to be. Thus these are examples of blow downs which do not rely upon positive scalar curvature to succeed. In addition, since the double branched cover bounds a negative definite manifold that is not a plumbing along trees, this example escapes the conclusion of \cite{Stip}, and thus it is open whether this family may give rise to symplectic blow down operations. \\

\begin{figure}
\includegraphics[width=4in,height=1.5in]{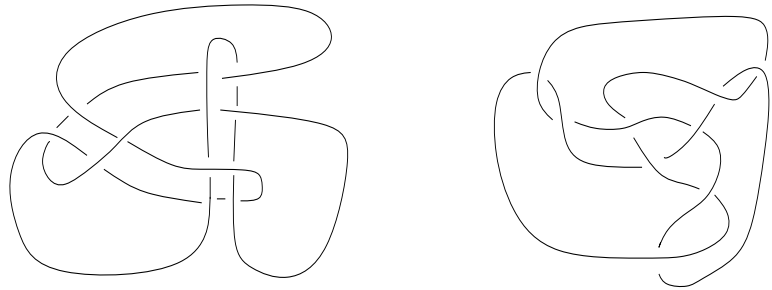}
\caption{The knot $9_{41}$ as a ribbon knot (left) and an alternating knot (right). This knot has an obvious $\Z/3\Z$ symmetry, but in fact it
has an additional $\Z/2\Z$-symmetry.}\label{fig:941}
\end{figure}

\begin{figure}
\includegraphics[width=4in,height=1.5in]{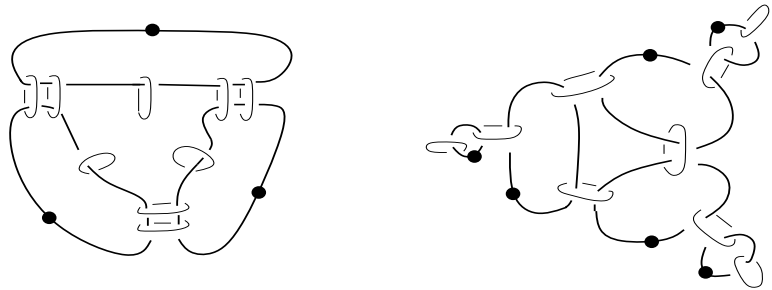}
\caption{The result of the construction from \cite{Doub} applied to $9_{41}$ (left) and $\overline{9}_{41}$ (right). All undotted unknots have framing $-1$.}\label{fig:dotted}
\end{figure}

\begin{figure}
\includegraphics[width=4in,height=1.5in]{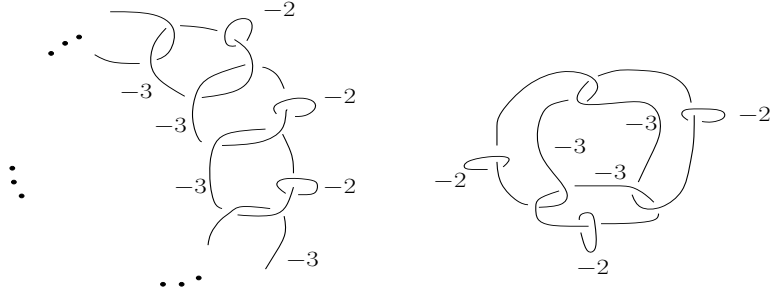}
\caption{The negative definite manifolds bounded by the branched double covers of knots in out family using the conventions of \cite{Unkn}. If the knot is the $n$-fold cylic cover, then there are $n$ two handles with framing $-3$, and $n$ two handles with framing $-2$. That for $\Si(9_{41})$ is shown on the right.}\label{fig:333222}
\end{figure}

\begin{figure}
\includegraphics[width=4in,height=1.5in]{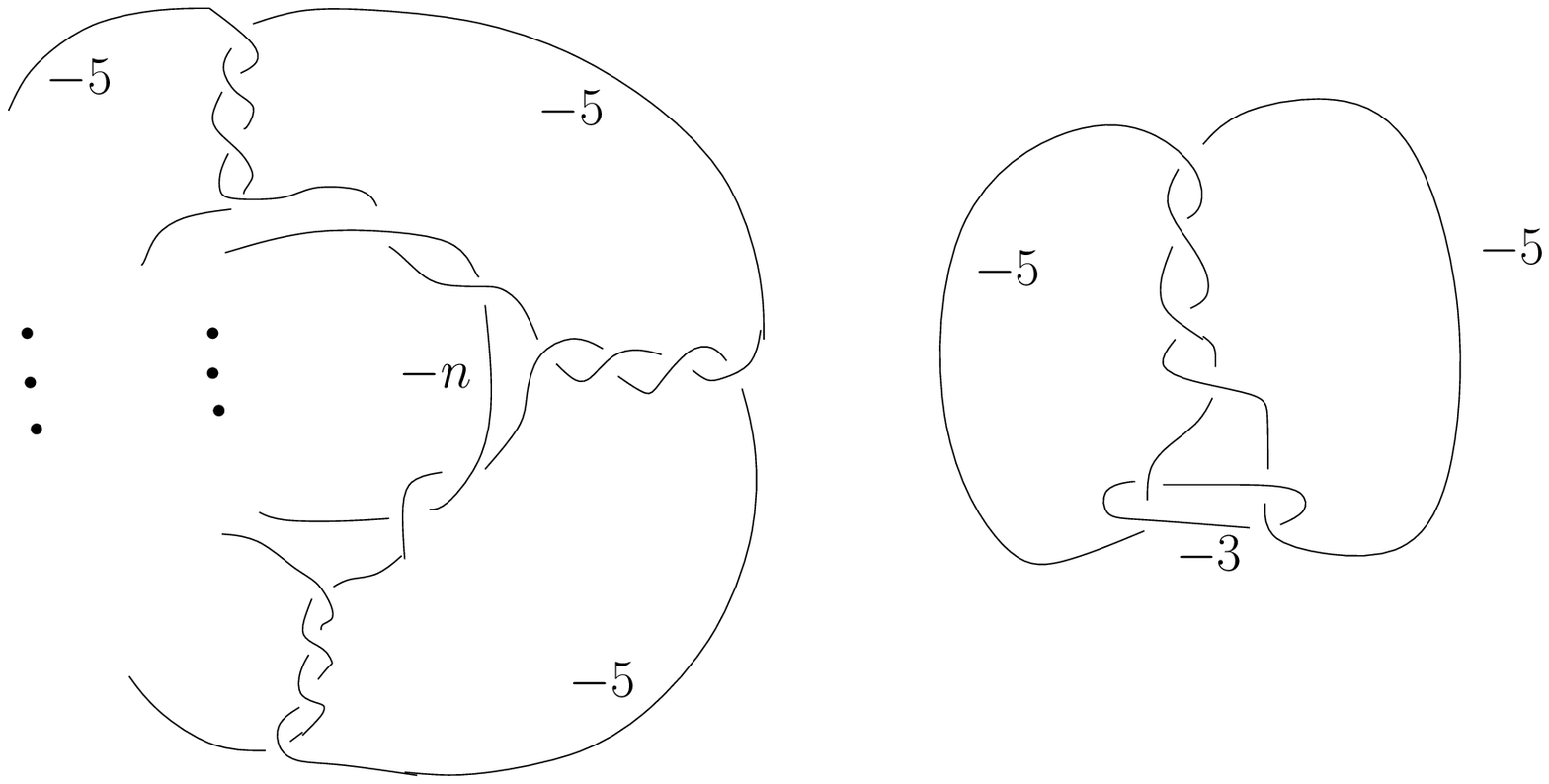}
\caption{The negative definite manifolds bounded by the branched double covers of the mirrors of knots in our family. If the knot is the $n$-fold cyclic cover, then there are $n-1$ handles with framing $-5$. That for $\Si(\overline{9}_{41})$ is shown on the right. }\label{fig:535}
\end{figure}
\ \\
\noindent To be specific, we depict the knot $9_{41}$ in Figure \ref{fig:941}. In \cite{Doub}, a construction is given for a sleek negative definite
manifold whose boundary is $\Si(9_{41})$. One forms the black/white coloring on the plane determined by the alternating projection and the convention that the over-strand at a crossing must cross a white region locally, when rotating counter-clockwise to get to the under-strand. We then associate one handles in our four manifold to the black regions by drawing a dotted circle in all but one black region, and lying in the plane. At each crossing
we add a $-1$ framed two handle to an unknot looping through the two dotted circles associated to the crossing. By judicious handleslides, we may arrange for a tree of dotted circles and unknots which can then be cancelled, leaving only two handles. That the result is negative definite is proven in \cite{Doub}. We show the initial stage of this construction for $9_{41}$ and $\overline{9}_{41}$ in Figure \ref{fig:dotted}. The results for the 
knots in our family are shown in Figures \ref{fig:333222} and \ref{fig:535}. \\
\ \\
\noindent For $9_{41}$ we check the vectors $K$ on the negative definite piece such that $|\langle K, v_{i} \rangle| \leq 3$ and $|\langle K, w_{i} \rangle | \leq 2$ ehre we have $v_{i} \cap w_{i}  = + 1$. The configuration is invariant with $S_{3}$ symmetry. The inverse of the intersection form
is then 
$$
I^{-1} = -\frac{1}{7}\left( \begin{array}{cccccc} 
		 6  & 4 & 4 & 3 & 2 & 2 \\
		 4  & 6 & 4 & 2 & 3 & 2 \\
		 4  & 4 & 6 & 2 & 2 & 3 \\
		 3  & 2 & 2 & 5 & 1 & 1 \\
		 2  & 3 & 2 & 1 & 5 & 1 \\
		 2  & 2 & 3 & 1 & 1 & 5 \\
		\end{array} \right) \\
$$
One then finds that $(-1,1,1,0,2,0)$ and its orbit all restrict as $Spin^{c}$ structures with $d(\mathfrak{s}) = 0$. Likewise for $(-1,-1,1,2,0,2)$ and
its orbit, and $(-1,1,1,0,2,0)$ restricts as the conjugate of the restriction of $(1,-1,-1,2,0,2)$. This last requires running through the equivalences
generated by the spheres in the diagram. Each vector gives rise to a path of length $6$ before arriving at the conjugate. The last set of vectors in this range which give $d(\mathfrak{s}) = 0$ are $(1,1,1,0,0,0)$, $(-1,-1,-1,0,0,0)$, $(1,1,1,-2,-2,-2)$ and $(-1,-1,-1,2,2,2)$, but these all restrict to the same $Spin^{C}$ structure. The equivalence in the last pair fits within the full path model, but the equivalence of either of the first two 
with the last two does not. Nevertheless, there is a class in the second homology whose Poincar\'e dual effects the equivalence (it is the sum of all
three of the $-3$ spheres), and this is all that is required in \cite{Doub}. With this caveat, we have found $13$ $Spin^{c}$ structures with the desired property. When we glue in the rational homology ball we obtain a $\Z/7\Z$ orbit which extends. There are two possibilities depending upon how we glue in the ball. These are found by using the $(1,1,1,-2,-2,-2)$ structure, which is totally invariant, as a basepoint since by conjugation it must be one of those to extend. Taking one of the other vectors, say  $(1,1,-1,0,2,0)$, we examine the subgroup generated by the difference with the basepoint. As $H_{1}(\Si(9_{41})) \cong \Z/7\Z\,\oplus\,\Z/7\Z$, this will be one of the possible subgroups. The result consists of the restrictions of
the following vectors, written in the order obtained:
$$
\begin{array}{c}
(1,1,-1,0,2,0) \ra (-1,1,1,0,0,2) \ra (-1,1,-1,0,2,2) \ra (1,-1,1,2,0,0) \\ \ra (1,-1,-1,2,2,0) \ra (-1,-1,1,2,0,2) \ra (1,1,1,0,0,0) \\
\end{array}
$$
\noindent The conjugation action reverses the order. The $\Z/3\Z$ rotational symmetry preserves this subgroup, taking the first element to the second to the fourth, and the third element to the fifth to the sixth. On the other hand, the $\Z/2\Z$ symmetry found from a $180^{\circ}$ rotation of the figure takes this subgroup into a subgroup containing the invariant structure and the other six vectors. Thus there are two possible sets of extensions
and which one depends upon the gluing through the $\Z/2\Z$ symmetry of $\Si(9_{41})$.   \\ 
\ \\
\noindent For $\overline{9}_{41}$, one has 

$$
I^{-1} = -\frac{1}{7}\left( \begin{array}{cccc} 
		 2 & 1  & 1 \\
		 1  & 3 & 1 \\
		 1  & 1  & 2\\
		\end{array} \right) \\
$$
\ \\
\noindent we can then compute the $K$-vectors with $|\langle K, v_{1} \rangle | \leq 5$, $|\langle K, v_{2} \rangle | \leq 5$, and $|\langle K, v_{3} \rangle | \leq 3$ which also give $d(\mathfrak{s}) = 0$. These are 

$$
\begin{array}{lcl}
\ &(-1,3,-1)\sim (1,-3,1) &\ \\
(1,-1,-3)& \hspace{.5in} &(-3,-1,1)\\
(-3,-1,3)& \hspace{.5in} &(3,-1,3) \\
(-3,1,-1)& \hspace{.5in} &(-1,1,-3) \\
(3,-1,1) &\hspace{.5in} &(1,-1,3) \\
(3,-1,-3)&\hspace{.5in}&(-3,1,2) \\
(-1,1,3) &\hspace{.5in}&(3,1,-1) \\
\end{array}
$$
\ \\
\noindent Since $H_{1}(\Si(\overline{9}_{41})) \cong \Z/7\Z\,\oplus\,\Z/7\Z$, we need to discover which set of $7$ will extend for a gluing
of the rational homology ball into the boundary. We can cut down the options a little by noticing that there are two symmetries, conjugation and
the obvious $\Z/2\Z$-symmetry. These imply that the spin structure extends; using it as a basepoint, we need only find $\Z/7\Z$ subgroups contained in the set of vectors with $d(\mathfrak{s}) = 0$. Any other vector in this set generates such a subgroup, and it is straightforward to check that the two columns are the two subgroups, taken one into the other by the symmetry. In fact, if we add $2(0,1,-2)$ to the spin structure, we can generate vectors equivalent to the first column, and these will specify the $Spin^{c}$ structures on the boundary by restriction. Note also that all of these vectors could occur for taut configurations. Using the handlebody diagrams, we can find Legendrian links producing Stein structures on the negative definite piece whose canonical classes sit within this set. \\
\ \\ 
\noindent {\bf IV.} There are two other slice, alternating, two bridge knots whose branched covers were not included by J. Park's results. These are $9_{27}$ and $10_{42}$ with double branched covers of $L(49,19)$ and $L(81,31)$ respectively. Neither are amphichiral, so their mirrors will also provide examples.The remaining alternating, slice knots that are also not two-bridge are $10_{87}$, $10_{99}$, and $10_{123}$. Of these, $10_{123}$ provides the next most likely example to show up in practice. Its branched double cover bounds the negative definite manifold found by $-3$ surgery on a ring of five unknots, each intersecting the next $+1$ times. However, we forgo the analysis for now. $10_{87}$ and $10_{99}$ have no such nice descriptions, but their double covers can also be found to bound rings of five unknots, with varying surgery coefficients.

\appendix

\section{Background on the Four-Manifold Invariants}\label{app:four}

\noindent We record here the general definition of the mixed invariants for a four dimensional cobordism. Let $W^{4}$ have boundary $-Y_{1} \cup Y_{2}$ and assume that $b_{\,2}^{+}(W) > 1$. Let $\mathfrak{s}$ be a $Spin^{c}$ structure on $W$. In \cite{Smoo} \Oz and \Sz define a map $F_{W,\mathfrak{s}}^{mix} : HF^{-}(Y_{1}, \mathfrak{s}_{1}) \ra HF^{+}(Y_{2}, \mathfrak{s}_{2})$. To do so, they decompose $W$ as $W_{1} \cup_{N} W_{2}$ where $b^{+}_{\,2}(W_{i}) > 0$ and $\delta H^{1}(N;\,\Z) = 0$ in $H^{2}(W, \partial W;\,\Z)$. Such decompositions are always available. In this setting one has the diagram: \\

$$
\begindc{\commdiag}[8]
\obj(0,0)[Ad]{$HF^{-}(Y_{1}, \mathfrak{s}_{1})$}
\obj(20,0)[Ad1]{$HF^{-}_{red}(N, \mathfrak{s}|_{N})$}
\obj(20,10)[Bl1]{$HF^{+}_{red}(N,\mathfrak{s}|_{N})$}
\obj(40,10)[Bl2]{$HF^{+}(Y_{2},\mathfrak{s}_{2})$}
\mor{Ad}{Ad1}{$F^{-}_{W_{1}, \mathfrak{s}_{W_{1}}}$}
\mor{Ad1}{Bl1}{$\tau_{N}^{-1}$}
\mor{Bl1}{Bl2}{$F^{+}_{W_{2},\mathfrak{s}|_{W_{2}}}$}
\enddc
$$
\ \\
\ \\
\noindent {\bf Note:} The conditions on $b_{\,2}^{+}$ ensure that the cobordism maps occur as presented. This map does not depend upon the particular choice of $N$ and defines a map equivariant with respect to the action of $\mathrm{Diff}^{+}(W)$ on $Spin^{c}$ structues. Furthermore, there are only finitely many $\mathfrak{s}$ for which this mixed map is not zero.\\
 \\
\noindent When $W$ occurs as the complement of two balls in a closed four manifold, \Oz and \Sz use this map to define their invariants of smooth four manifolds, $\Phi_{X,\mathfrak{s}}$, by

$$
  F^{mix}_{X - (B_{1} \cup  B_{2}), \mathfrak{s}}((U^{n} \otimes \gamma) \cdot\Theta^{-}_{(-2)}) = \pm \Phi_{X,\mathfrak{s}}(U^{n}\otimes \gamma) \Theta^{+}_{(0)} + (\mathit{independent\ terms})  
$$

\noindent where $\Theta_{(\cdot)}$ is the standard generator of the homology for $S^{3}$ in that grading, within the $\pm$-homology as appropriate. Here $\gamma$ is an element in $\Lambda^{\ast}(H_{1}(X)/\mathrm{Tors})$ and the map is zero except on those elements whose degree equals:

$$
D(X, \mathfrak{u}) = \frac{c_{1}(\mathfrak{u})^{2} - (2\,\chi(X) +  3\,\sigma(X))}{4}
$$

\noindent where $U$ has degree $2$ and an element of $H_{1}(X;\,\Z)$ has degree $1$.\\
  \\
\noindent We will use the map when $Y_{1}$ is not $S^{3}$, but $Y_{2}= S^{3}$. We can then define

$$
  F^{mix}_{X, \mathfrak{s}} ((U^{n} \otimes \gamma) \cdot {\bf X})= \pm \Phi_{X,\mathfrak{s}}(U^{n}\otimes \gamma \otimes {\bf X})\Theta_{(0)}^{+} + (\mathit{independent\ terms})
$$

\noindent where ${\bf X}$ is an element of $HF^{-}(Y_{1}, \mathfrak{s}_{1})$. \\
 \\
\noindent If $Y_{1}$ is a rational homology sphere and $W_{0}$ has boundary equal to $Y_{1}$ with $b_{1}(W) = 0$ then the formula for compositions of cobordism maps, \cite{Smoo}, implies that

$$
\Phi_{W \cup X, \mathfrak{s} \# \mathfrak{s}_{0}}(U^{n + m}\otimes \gamma) = \Phi_{X, \mathfrak{s}}( U^{n}\otimes \gamma \otimes F^{-}_{W, \mathfrak{s}_{0}}(U^{m} \cdot \Theta^{-}_{(-2)}))
$$

\pagebreak

\begin{figure}
\includegraphics[height=7in,width=5in]{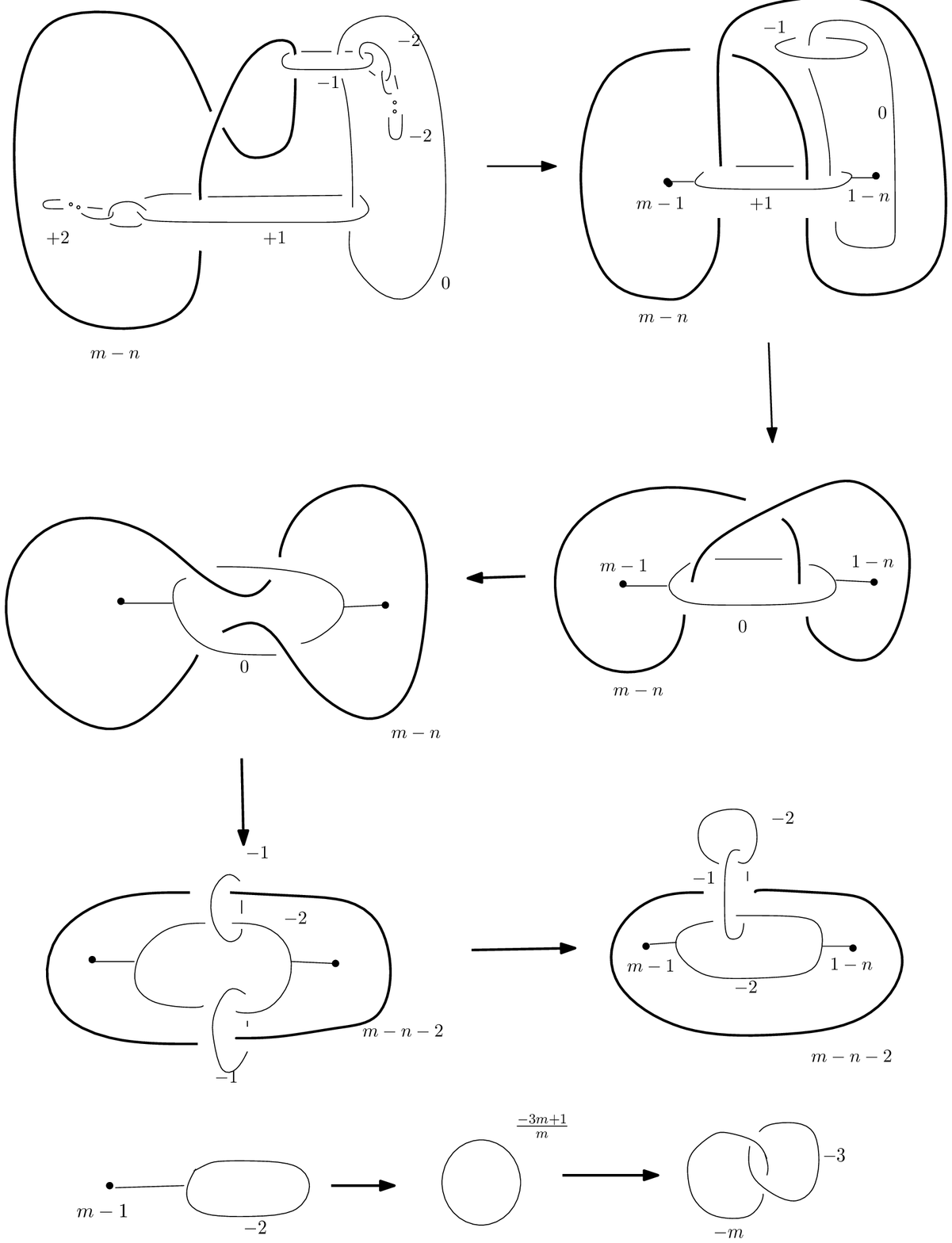}
\caption{The steps in showing the boundary equivalence of Figure \ref{fig:257fam} and Figure \ref{fig:257neg}. A dot with a number indicates a
linear chain of $\pm 2$ spheres. If the number is positive there are that many $+2$ spheres; if the number is negative there are minus that number of $-2$ spheres. The sequence at the bottom shows how to ultimately remove the $+2$ spheres in favor of negative self-intersections.} 
\end{figure}


\begin{thebibliography}{99}
\bibitem{Cass}A. J. Casson \& J. L. Harer, {\em Some Homology Lens Spaces Which  Bound Rational Homology Balls}. Pac. J. of Math. 96(1): 23-36, 1981.

\bibitem{Fint}R. Fintushel \& R. Stern, {\em Rational Blow-downs of Smooth Four Manifolds}. J. Diff. Geom. 46:181-235, 1997.

\bibitem{Froy}K. Froyshov, {\em A Generalized Blow-up Formula for Seiberg-Witten Invariants}. math.DG/0604242 v2.

\bibitem{Four}R. E. Gompf \& Andras I. Stipsicz, {\underline{4-Manifolds and Kirby Calculus}}. Graduate Studies in Mathematics v. 20. AMS. Providence, Rhode Island, 1999.

\bibitem{Kirb}J. Harer, A. Kas, \& R. Kirby, {\em Handlebody Decompositions of Complex Surfaces}. Memoirs of the AMS, \#350. AMS. Providence, Rhode Island, July 1986.

\bibitem{Mark}S. Jabuka \& T. Mark, {\em Gluing Results in Heegaard-Floer Homology}, preprint.

\bibitem{Mich}M. Michalogiorgaki, {\em Rational Blow-down Along Wahl Type Plumbing Trees of Spheres and Exotic Smooth Structures on $CP^2\#k\bar{CP^2},k \in {6,7,8,9}$}. math.GT/0607608. 

\bibitem{Naik}S. Naik, {\em Equivariant Concordance of Knots in $S^{3}$}. Proceedings of Knots '96, Ed. Shin'ichi Suzuki. World Scientific Publishing
Co. 1997, 81-89.

\bibitem{Neum}W. D. Neumann, {\em An Invariant of Plumbed Homology Spheres}. Topology Symposium, Siegen 1979. LNM \#788: 125-144.

\bibitem{Ngky}K. Y. Ng, {\em Groups of Ribbon Knots}

\bibitem{Nico}L. I. Nicolaescu, {\em Notes on Seiberg-Witten Theory}. AMS, Providence, Rhode Island, 2000.

\bibitem{Owen}Owens \& Strle, {\em Rational Homology Spheres and the Four Ball Genus}. math.GT/0308073.

\bibitem{AbsG}P. Ozsv{\'a}th \& Z. Szab{\'o}, {\em Absolute Grading Floer Homologies and Intersection Forms for Four-Manifolds with Boundary}. Advances in Mathematics, 173 (2): 179-261, 2003.

\bibitem{Symp}P. Ozsv{\'a}th \& Z. Szab{\'o}, {\em Holomorphic Triangle Invariants and the Topology of Symplectic Four-Manifolds}. Duke Math. J. 121(1), 1-34, 2004.

\bibitem{Hom3}P. Ozsv{\'a}th \& Z. Szab{\'o}, {\em Holomorphic Disks and Topological Invariants for Closed Three Manifolds}.  Ann. of Math. (2)  159(3):1027-1158, 2004.

\bibitem{3Man}P. Ozsv{\'a}th \& Z. Szab{\'o}, {\em Holomorphic Disks and Three Manifold Invariants: Properties and Applications}. Ann. of Math. (2) 159(3): 1159-1245, 2004.

\bibitem{Smoo}P. Ozsv{\'a}th \& Z. Szab{\'o}, {\em Holomorphic Triangles and Invariants of Smooth Four Manifolds}. math.SG/0110169.

\bibitem{Unkn}P. Ozsv{\'a}th \& Z. Szab{\'o}, {\em Knots with Unknotting Number One and Heegaard Floer Homology}. math.GT/0401426 v1.

\bibitem{Lens}P. Ozsv{\'a}th \& Z. Szab{\'o}, {\em On Knot Floer Homology and Lens Space Surgeries}. Topology 44(6):1281-1300, 2005.

\bibitem{Plum}P. Ozsv{\'a}th \& Z. Szab{\'o}, {\em On the Floer Homology of Plumbed Three-Manifolds}. Geom. \& Topol. 7:185-224, 2003.

\bibitem{Doub}P. Ozsv{\'a}th \& Z. Szab{\'o}, {\em On the Heegaard Floer Homology of Branched Double Covers}. math.GT/0309170 v1.

\bibitem{Park}J. Park, {\em Seiberg-Witten Invariants of Generalised Rational  Blow-Downs}. Bull. Austal. Math. Soc. V.56: 363-384, 1997.

\bibitem{Rust}R. Rustamov, {\em On the Heegaard-Floer Homology of Plumbed Three Manifolds with $b_{1} = 1$}. math.SG/04051181 v1.

\bibitem{Sala}D. Salamon, {\em Spin Geometry and Seiberg-Witten Invariants}, preprint, 1996.

\bibitem{Stip}A. Stipsicz \& Z. Szab{\'o}, {\em A Note on Symplectic Rational Blow Downs}. math.GT/0511101 v1.
 
\bibitem{Wahl}J. Wahl, {\em Smoothings of Normal Surface Singularities}. Topology v. 20:219-246, 1981.

\end{thebibliography}
\end{document}